\documentclass[a4paper,12pt]{article}
\usepackage{latexsym}
\usepackage{amsmath}
\usepackage{amsthm}

\usepackage{txfonts}
\usepackage{bm}
\usepackage{color}	
\usepackage{epic,eepic}

\newcommand{\REV}[1]{{\color{red}#1}} 
\newcommand{\RED}[1]{{\color{red}#1}} 
\newcommand{\BLU}[1]{{\color{blue}#1}}

\renewcommand{\REV}[1]{{\color{black}#1}} 
\renewcommand{\RED}[1]{{\color{black}#1}} 
\renewcommand{\BLU}[1]{{\color{black}#1}}

\usepackage{geometry}
\numberwithin{equation}{section}

\newtheorem{theorem}{Theorem}[section]
\newtheorem{proposition}[theorem]{Proposition}
\newtheorem{corollary}[theorem]{Corollary}
\newtheorem{lemma}[theorem]{Lemma}
\newtheorem{claim}[theorem]{Claim}

\newtheorem{remark}{Remark}[section]

\newcommand{\memo}[1]{{\bf \small \RED{[\bf MEMO:}} \BLU{ #1}  {\bf \small \RED{end]} }}  
   \renewcommand{\memo}[1]{}           
\newcommand{\OMIT}[1]{{\bf [OMIT:} #1 \ {\bf --- end OMIT] }}  
   \renewcommand{\OMIT}[1]{}            

\newcommand{\finbox}{\hspace*{\fill}$\rule{0.17cm}{0.17cm}$}

\newcommand{\odotZ}{\overset{....}}

\newcommand{\BB}{\hspace*{\fill}$\rule{0.17cm}{0.17cm}\ \rule{0.17cm}{0.17cm}$}

\newcommand{\Proof}{\noindent {\bf Proof.  }}

\newcommand{\llceil}{\bigg\lceil} 
\newcommand{\rrceil}{\bigg\rceil} 

\newcommand{\llfloor}{\bigg\lfloor} 
\newcommand{\rrfloor}{\bigg\rfloor} 

\usepackage{lineno}
\newcommand*\patchAmsMathEnvironmentForLineno[1]{
  \expandafter\let\csname old#1\expandafter\endcsname\csname #1\endcsname
  \expandafter\let\csname oldend#1\expandafter\endcsname\csname end#1\endcsname
  \renewenvironment{#1}
     {\linenomath\csname old#1\endcsname}
     {\csname oldend#1\endcsname\endlinenomath}}
\newcommand*\patchBothAmsMathEnvironmentsForLineno[1]{
  \patchAmsMathEnvironmentForLineno{#1}
  \patchAmsMathEnvironmentForLineno{#1*}}
\AtBeginDocument{
\patchBothAmsMathEnvironmentsForLineno{equation}
\patchBothAmsMathEnvironmentsForLineno{align}
\patchBothAmsMathEnvironmentsForLineno{flalign}
\patchBothAmsMathEnvironmentsForLineno{alignat}
\patchBothAmsMathEnvironmentsForLineno{gather}
\patchBothAmsMathEnvironmentsForLineno{multline}
}

\begin{document}

\title{A Discrete Convex Min-Max Formula  \\  for Box-TDI Polyhedra}

\author{Andr\'as Frank%
\thanks{MTA-ELTE Egerv\'ary Research Group,
Department of Operations Research, E\"otv\"os University, P\'azm\'any
P.~s.~1/c, Budapest, Hungary, H-1117. 
e-mail:  {\tt frank@cs.elte.hu}. 
ORCID: 0000-0001-6161-4848.
}
 \ \ and \ 
{Kazuo Murota%
\thanks{Department of Economics and Business Administration,
Tokyo Metropolitan University, Tokyo 192-0397, Japan, 
e-mail:  {\tt murota@tmu.ac.jp}. 
ORCID: 0000-0003-1518-9152.
}}}


\date{July 2020 / August 2020 / January 2021}

\maketitle

\begin{abstract}
A min-max formula is proved for the minimum of an
integer-valued separable discrete convex function where the minimum is
taken over the set of integral elements of a box total dual integral (box-TDI) polyhedron.  
One variant of the theorem uses the notion of conjugate function 
(a fundamental concept in non-linear optimization)
but we also provide another version that avoids conjugates, and its
spirit is conceptually closer to the standard form of classic min-max
theorems in combinatorial optimization.  
The presented framework provides a unified background for separable convex minimization over
the set of integral elements of the intersection of two integral
base-polyhedra, submodular flows, L-convex sets, and polyhedra defined
by totally unimodular (TU) matrices.  
As an unexpected  application, we show how a wide class of  inverse
combinatorial optimization problems can be covered by this new framework.
\end{abstract}





\newpage

\tableofcontents

\newpage



\section{Introduction}
\label{intro}

A central aspect of convex optimization is minimizing a convex function over a convex set.  
Discrete convex analysis \cite{Murota98a,Murota03} considers discrete convex functions.  
It turned out that there are two strongly interrelated general classes,
M-convex and L-convex functions, for which fundamental min-max
theorems can be formulated.  
It is important to distinguish between the cases 
when we minimize over real or over integer vectors.  
For example, one may be interested in finding 
a minimum 
\RED{
$l_{2}$-norm
}
element of an integral base-polyhedron $B$ (say) 
or a minimum 
\RED{
$l_{2}$-norm
}
{\rm integral} element of $B$.  
These are pretty different problems as the continuous version has a unique solution \cite{Fujishige80},
while 
the set of integral optima \cite{FM19partA} concerning base-polyhedra has a rich structure.  
In the present work, we discuss the second type of minimization 
when the function to be minimized is an integer-valued
separable discrete convex function.  
It was proved in \cite{Murota03}
that these functions are exactly those which are both 
M$\sp{\natural}$-convex and L$\sp{\natural}$-convex.  
In this sense separable discrete convex functions are rather special 
but this speciality makes it possible that we can develop min-max theorems  when 
we minimize over a discrete box-TDI set, 
a much wider class than M$\sp{\natural}$-convex or L$\sp{\natural}$-convex sets.  
Box-TDI linear systems and polyhedra 
(defined formally below) were introduced by Edmonds an Giles \cite{Edmonds-Giles84}, 
studied in detail by Cook \cite{Cook83,Cook86}, 
and recently by Chervet, Grappe, and Robert \cite{CGR}.  
We shall call the set of integral elements of an integral box-TDI
polyhedron a {\bf discrete box-TDI set}, or just a {\bf box-TDI set}.

Our main goal is to develop a general min-max formula for the minimum
of an integer-valued separable discrete convex function $\Phi$ over a
discrete box-TDI set.  Actually, we exhibit two equivalent forms.  
One of them makes use of the discrete version of Fenchel conjugate, 
a fundamental concept from non-linear (continuous) optimization 
(see \cite{BL,HL01,Roc70}).  
But we also develop another form which does not rely on 
the concept of conjugate, and therefore this version is
conceptually closer to classic min-max theorems of combinatorial optimization 
like the ones of Menger, K{\H o}nig, Egerv\'ary,
Dilworth, Ford$+$Fulkerson, Tutte, Edmonds, Lucchesi$+$Younger, etc.

Our general framework includes as a special case the corresponding
optimization problems for totally unimodular (TU) matrices, 
in particular, circulations and 
\REV{
tensions (= potential-differences).  
}
The results can also be applied to submodular flows, 
in particular to the intersection of two base-polyhedra.  
As a special case, we derive a min-max theorem for
the minimum square-sum of an integer-valued (!) feasible circulation or maximum flow.

It is our important goal to bring those readers closer to 
discrete convex optimization who are not particularly familiar with the notion of conjugate.  
The present work, apart from one exception, does not deal with algorithmic issues, 
but we hope that our min-max formulas
pave the way to forthcoming researches for constructing strongly
polynomial algorithms to compute the optima in question.

As an unexpected application, we shall show in Section \ref{inverz}
how a significant part of inverse combinatorial optimization problems
can be modelled in this new framework.  
We provide a min-max theorem for the minimum total change (measured in 
\RED{
$l_{1}$-norm)
}
of a given cost-function $w_{0}$ for which a specified element of a discrete box-TDI set
(for example, a spanning tree of a graph) 
becomes a cheapest one with respect to the modified cost-function $w$.  
Even the more general inverse problem fits into our framework 
when each element from a specified list is expected 
to be a cheapest one with respect to the 
\REV{
desired
}
$w$.

In the present work, for the sake of technical simplicity, we
concentrate on integer-valued functions.  
It should, however, be emphasized that all the results can be 
extended in a natural way to real-valued separable discrete convex functions, as well.

\subsection{Notions and notation}

Let ${\bf R}$, ${\bf Q}$, and ${\bf Z}$ denote the set of 
reals, rationals, and integers, respectively.  
When it does not make any confusion, 
we do not distinguish between row- and column-vectors.  
For example, if $u$ and $v$ are vectors from ${\bf R}\sp{n}$, then $uv=vu$
denotes their scalar product.  For a vector $w$, we use the notation
$w\sp{2}$ for the scalar product $ww$, and will refer to $w\sp{2}$ as
the {\bf square-sum} of $w$.  If $Q$ is an $m$-by-$n$ matrix while
$x\in {\bf R}\sp{n}$ and $y\in {\bf R}\sp{m}$ are vectors, then $x$ is
considered a column-vector in the product $Qx$, 
while $y$ is considered a row-vector in $yQ$.

Throughout we work with a ground-set $S$ with $n$ elements.  
The incidence or characteristic vector 
\REV{
of a subset $X$ of $S$ is denoted by $\chi_{X}$, and 
$\chi_{S}$ will be briefly denoted by ${\bf \underline{1}}$. 
} 
For elements $s,t\in S$, we call a subset
$X\subset S$ an {\bf $s\overline{t}$-set} if $s\in X\subset S-t$.  For a
function $f$ on $S$, the set-function $\widetilde f$ is defined by
$\widetilde f(X):=\sum [f(s):s\in X]$ $(X\subseteq S)$.

For a polyhedron 
$R:=\{x:  Qx\geq p\} \subseteq {\bf R}\sp{S}$, 
$\odotZ{R}$ denotes the set of integral elements of $R$, that is,
\begin{equation} \label{odottR}
\odotZ{R}:= R \cap {\bf Z} \sp{S}. 
\end{equation}
\noindent 
For a cost-function $w$ on $S$, let $\mu_{R}(w)$ denote the minimum of 
$\{wx:  x\in R\}$, while 
\REV{
$\mu_{\odotZ{R}}(w) := \min \{wx:  x\in \odotZ{R}\}$. 
}
We say that an element $z\sp{*}$ of $R$ is a 
{\bf $w$-minimizer} 
if $wz\sp{*} \leq wx$ holds for every $x\in R$, that is, $wz\sp{*}=\mu_{R}(w)$.

The {\bf effective domain} \cite{Murota03,Roc70} (or sometimes just
{\bf domain} \cite{BL,HL01}) ${\rm dom}(\varphi )$ 
of an integer-valued function 
$\varphi :{\bf Z} \rightarrow {\bf Z} \cup \{-\infty ,+\infty \}$ 
is the set of integers where $\varphi$ is finite.  
When we say that a function $\varphi$ is integer-valued, 
we allow that some values of $\varphi$ may be $-\infty $ or $+\infty $. 
A function
$\varphi :{\bf Z} \rightarrow {\bf Z} \cup \{+\infty \}$ 
is called {\bf discrete convex} if 
\begin{equation} \label{univerdcfndef}
 \varphi (k-1) + \varphi (k+1) \geq 2\varphi (k)
\end{equation}
for each $k\in {\rm dom} (\varphi)$.  
Let $\varphi '$ denote the function defined on ${\bf Z}$ by
\begin{equation} 
\varphi '(k):= \varphi (k+1)-\varphi (k) \qquad (k\in {\bf Z}).
\label{(fiprime)} 
\end{equation}
\noindent 
The function $\varphi '$ may intuitively be considered 
the discrete right derivative of $\varphi$. 
Clearly, $\varphi$ is discrete convex precisely 
if $\varphi '$ is monotone non-decreasing.
The effective domain of a discrete convex function is the set of integers 
in a (possibly unbounded) interval.

When we are given a function $\varphi_{s}$ for every $s\in S$, 
the functions 
$\Phi:  {\bf Z}\sp{S} \rightarrow {\bf Z} \cup \{+\infty \}$
and 
$\Phi':  {\bf Z}\sp{S} \rightarrow {\bf Z} \cup \{-\infty ,+\infty \}$ 
are defined by:
\begin{equation} 
\Phi(z):= \sum_{s\in S} \varphi_{s}(z(s)), \qquad
\Phi'(z):= \sum_{s\in S} \varphi'_{s}(z(s)).  
\label{(Phidef)} 
\end{equation}
\noindent 
When each $\varphi_{s}$ is discrete convex, 
$\Phi$ is called a {\bf separable discrete convex} function.  
The {\bf discrete conjugate} function 
$\varphi \sp{\bullet}$ of a function 
$\varphi :{\bf Z} \rightarrow {\bf Z} \cup \{+\infty \}$ 
is defined for any integer $\ell$ by 
\begin{equation} 
\varphi \sp{\bullet}(\ell) := \max \{k\ell- \varphi (k) :  k\in {\bf Z}\}, 
\label{(conj-def)} 
\end{equation} 
while the discrete conjugate
$\Phi \sp{\bullet}$ of $\Phi$ is defined for $w\in {\bf Z} \sp{S}$ by
\[
 \Phi \sp{\bullet}(w):= \sum_{s\in S} \varphi \sp{\bullet}_{s}(w(s)).
\]
\noindent 
Note that 
$\Phi \sp{\bullet}(w) = \max \{wz - \Phi(z):  z\in {\bf Z} \sp{S} \}$, 
and this latter expression is actually the definition of the
discrete conjugate of an arbitrary integer-valued function $\Phi$ on ${\bf Z}\sp{S}$.

Note that $\varphi \sp{\bullet}(\ell)$ may be $+\infty $ 
(when $\{k\ell- \varphi (k) :  k\in {\bf Z}\}$ is not bounded from above) 
and hence using supremum would be formally a bit more precise but we keep the term maximum.  
It should be emphasized that in the original definition
of Fenchel conjugate in continuous optimization \cite{BL}, the maximum
is taken over all real values $k$ and not only on integer $k$'s.

Let $p:2\sp{S}\rightarrow {\bf Z} \cup \{-\infty \}$ 
be an integer-valued (fully) supermodular function on a ground-set $S$ 
for which the value $p(S)$ is finite.  
When we say that a function $p$ is supermodular, 
we always mean that the supermodular inequality
$p(X)+p(Y)\leq p(X\cap Y)+p(X\cup Y)$ 
holds for every pair $\{ X,Y \}$ of subsets of $S$.  
Since weaker supermodular functions (e.g., intersecting, crossing) 
are also important in applications, sometimes
we (over-) emphasize by saying that $p$ is \lq fully\rq \ supermodular.

Let 
\[
  B:=B'(p):= \{x:  \widetilde x(Z)\geq p(Z) \ \mbox{ for every } \  Z\subset S, 
 \mbox{ and } \  \widetilde x(S)=p(S)\} 
\]
be the {\bf base-polyhedron} defined by $p$.  
\REV{
Since $p$ is integer-valued, $B$ is an integral polyhedron,
which, in turn, determines $p$ uniquely as $p(Z)=\min\{ \widetilde x(Z):  x \in B \}$.
}
Note that the complementary function $\overline{p}$, 
defined by $\overline{p}(X):=p(S)-p(S-X)$, 
is submodular and 
$B'(p) = B(\overline{p}) :=\{x:  \widetilde x(Z)\leq \overline{p}(Z)$ for every 
$Z\subset S$, and $\widetilde x(S)=\overline{p}(S)\}$.  
That is, a base polyhedron can be defined by a submodular function as well.

In discrete convex analysis \cite{Murota03}, 
the set $\odotZ{B}$ of integral elements of $B$ is called an {\bf M-convex set} 
and the intersection of two M-convex sets an {\bf M$_{2}$-convex set}.  
A fundamental theorem of Edmonds \cite{Edmonds70} states that a set is
M$_{2}$-convex precisely if it is the set of integral elements of the
intersection of two integral base-polyhedra.

\subsection{Starting points}

A starting point of the present work is the problem of
finding/characterizing an element of an M-convex set $\odotZ{B}$ for
which an integer-valued separable discrete convex function $\Phi(z)$
in \eqref{(Phidef)} is minimum.  
\REV{
It is a basic property of integral base-polyhedra 
(see, e.g.,  \cite{Frank-book}) 
that the intersection of an integral box with an integral base-polyhedron is itself 
an integral base-polyhedron.  
}
Since the effective domain of $\Phi$ is a box, 
it follows that we can replace $B$ with this intersection, or in other words, 
we may assume that $\Phi$ is finite-valued on the whole M-convex set $\odotZ{B}$.

A min-max theorem for separable discrete convex functions on an
M-convex set can be obtained as a special case of the Fenchel-type
discrete duality theorem \cite{Murota03} (Theorem 8.21) concerning
discrete convex functions which are not necessarily separable.  
The formulation needs the well-known concept of linear (or Lov\'asz)
extension $\hat p$ of $p$ which is recalled in 
\REV{
\eqref{lovextdef} in Section~\ref{SCmsetm2set}.  
}
We also hasten to recall a basic theorem of Edmonds \cite{Edmonds70,Edmonds71} 
asserting that 
$\hat p(w)= \min \{ wz:  z\in \odotZ{B}\} \ (= \min \{wz:  z\in B\})$.  
For an element $z\in \odotZ{B}$, 
we call a subset $X\subseteq S$ {\bf $z$-tight} if
$\widetilde z(X) = p(X)$.  For a vector $w\in {\bf Z}\sp{S}$, 
we call a non-empty set $X\subseteq S$ a {\bf strict $w$-top set} 
if $w(s)>w(t)$ holds whenever $s\in X$ and $t\in S-X$.  
Note that the strict $w$-top sets form a chain.

Recall that an M-convex set $\odotZ{B}$ was defined as the set of
integral elements of an integral base-polyhedron $B$, that is,
\begin{equation} 
\odotZ{B} := B\cap {\bf Z} \sp{S}. 
\label{(odottB)} 
\end{equation}
\REV{
Although the present work was highly motivated by the theory of
discrete convex analysis (DCA) \cite{Murota03}, especially in
formulating some of the theorems, we do not rely on any prerequisite
from DCA, that is, each of our results and proofs are direct and
self-contained.  For DCA experts, however, as well as for readers who
may want to get acquainted with DCA in the future, it may be
beneficial if we point out some links to DCA.  For example, the
following three theorems were originally proved with tools from DCA.
It will be one of our goals to derive them directly (in a more general form).  
}

\begin{theorem} [\cite{Frank-Murota.2}] \label{min-max.orig} 
Suppose that an integer-valued separable discrete convex function $\Phi$ 
is finite-valued and bounded from below on an M-convex set $\odotZ{B}$
defined by an integer-valued (fully) supermodular function $p$
(allowing $-\infty $ values).  Then
\begin{equation} 
\min \{ \Phi(z) :  z\in \odotZ{B} \} = 
 \max \{\hat p(w) - \Phi \sp{\bullet}(w) :  \ w\in {\bf Z}\sp{S} \}, 
\label{(min-max.orig)} 
\end{equation}
where $\Phi\sp{\bullet}$ denotes the discrete conjugate of $\Phi$ and
$\hat p$ denotes the linear extension of $p$ 
(and hence $\hat p(w)=\mu_B(w)$).  
Moreover, an element $z\sp{*}\in \odotZ{B}$ is a
$\Phi$-minimizer if and only if there is an integer-valued function
$w\sp{*}$ on $S$ meeting the following optimality criteria:
\begin{align} 
& 
\mbox{\rm each strict $w\sp{*}$-top set is $z\sp{*}$-tight, }\
\label{(optcrit1b)} 
\\
& \mbox{\rm $ \varphi_{s}'(z\sp{*}(s)-1) \leq w\sp{*}(s) \leq \varphi_{s}'(z\sp{*}(s))$ \ for each $s\in S$,}\ 
\label{(optcrit2d)} 
\end{align}
\noindent 
or writing \eqref{(optcrit2d)} concisely:
\begin{equation} 
\Phi '(z\sp{*} - {\bf \underline{1}} ) \   \leq \   w\sp{*}   \ \leq \   \Phi' (z\sp{*}).  
\label{(optcrit2.conc} 
\end{equation} 
\finbox
\end{theorem}

Actually, the general Fenchel-type min-max theorem in \cite{Murota03}
also implies the following extension of 
Theorem~\ref{min-max.orig} to M$_{2}$-convex sets.

\begin{theorem} [\cite{Frank-Murota.2}] \label{min-max-M2} 
Let $B_{1}:=B'(p_{1})$
and $B_{2}:=B'(p_{2})$ be base-polyhedra defined 
by integer-valued supermodular functions $p_{1}$ and $p_{2}$ 
for which $B:=B_{1}\cap B_{2}$ is non-empty.  
Let $\Phi$ be a finite integer-valued separable discrete convex function on $B$ 
which is bounded from below on $B$.  
Then one has:
\begin{equation} 
\min \{\Phi(z) :  z\in \odotZ{B} \} 
= \max \{\hat p_{1}(w_{1})
+ \hat p_{2}(w_{2}) - \Phi\sp{\bullet}(w_{1}+w_{2}) :  \ w_{1}, w_{2}\in {\bf Z}\sp{S}\}. 
\label{(min-max.M2)} 
\end{equation} 
\finbox
\end{theorem}

In Section \ref{spec.poly}, we shall derive these theorems from the
new min-max formula concerning discrete box-TDI sets.  
It is worth mentioning already at this point that in important special cases 
the discrete conjugate of $\Phi$ can be explicitly given.  
For example, let
$\Phi(z):= z\sp{2}$ \ ($= \sum [z(s)\sp{2} :  s\in S]$).
For any real number $\alpha \in {\bf R}$, 
let $\lfloor \alpha \rfloor $ denote the largest integer not larger than $\alpha$, 
and $\lceil \alpha \rceil$ the smallest integer not smaller than $\alpha$. 
Then Theorems \ref{min-max.orig} and \ref{min-max-M2} can be specialized, as follows.

\begin{theorem}[\cite{Frank-Murota.2,FM19partA}]
\label{minsq-M2} 
Let $B=B'(p)$ be an integral base-polyhedron.  Then
\begin{equation} 
\min \{ z\sp{2} :  z\in \odotZ{B} \} 
\ = \ 
\max \{\hat p(w) - \sum_{s\in S} 
  \llfloor  \frac{w(s)}{ 2 }\rrfloor  
  \  \llceil \frac{w(s)}{ 2 } \rrceil
:  \ w \in {\bf Z}\sp{S}\}.  
\label{(minsq.M)} 
\end{equation}
\noindent 
Let $B_{1}:=B'(p_{1})$ and $B_{2}:=B'(p_{2})$ be integral base-polyhedra 
for which $B:=B_{1}\cap B_{2}$ is non-empty.  Then
\begin{align} 
& \min \{ z\sp{2} :  z\in \odotZ{B} \} 
\nonumber \\ & =
\max \{\hat p_{1}(w_{1}) + \hat p_{2}(w_{2}) - \sum_{s\in S} 
  \llfloor  \frac{w_{1}(s) + w_{2}(s)}{ 2 }\rrfloor  
  \  \llceil \frac{ w_{1}(s) + w_{2}(s)}{ 2 }\rrceil 
 :  \ w_{1}, w_{2}\in {\bf Z}\sp{S}\}.  
\label{(minsq.M2)} 
\end{align} 
\finbox
\end{theorem}

These results were formulated first in \cite{Frank-Murota.2} with a
proof relying on the general discrete Fenchel-type duality theorem \cite{Murota03}.  
We shall directly derive not only Theorems \ref{min-max.orig} and \ref{min-max-M2} 
but a variant, as well, which does not use the concept of conjugate.  
Furthermore, we shall show that the role of the M-convex or M$_{2}$-convex set 
in these theorems is only that they are discrete box-TDI sets.  
Note that it is a basic property of base-polyhedra that they are box-TDI and a theorem of
Edmonds and Giles \cite{Edmonds-Giles} implies that the intersection
of two base-polyhedra is also a box-TDI polyhedron.  
Therefore our main min-max theorem concerning discrete box-TDI sets will 
imply these special cases.

As mentioned above, the present work does not consider algorithmic
aspects, apart from one exception.  
In Section \ref{min-max.proof}, 
we shall provide an algorithmic approach to compute 
the dual optimum in Theorem~\ref{min-max.orig}, 
but even that algorithm can work only if a primal optimal solution is already available.  
But constructing a strongly polynomial algorithm for computing the primal optimum (that
is, a $\Phi$-minimizer element of an M-convex set) already in the
special case of weighted 
\RED{
square-sum 
}
(when $\Phi(z):= \sum [c(s)w(s)\sp{2}:s\in S]$, each $c(s)$ is positive) 
remains a major research problem.  
In the more general Theorem~\ref{min-max-M2}, 
the even more special case when $\Phi(w)= w\sp{2}$ is wide open from an
algorithmic point of view.


\section{Box-TDI systems and polyhedra}
\label{boxTDIsyspoly}

In what follows, $Q$ is an integral matrix and $p$ is an integral vector.  
Throughout we assume that there is a one-to-one correspondence 
between the columns of $Q$ and the elements of ground-set $S$.

Edmonds and Giles \cite{Edmonds-Giles,Edmonds-Giles84} called a
(rational) linear system $Qx\geq p$ {\bf totally dual integral} 
({\bf TDI}) if the maximum in the linear programming duality equation
\begin{equation}  
\min\{cx:  Qx\geq p\} \ = \ \max \{yp:  y\geq 0, yQ=c \}
\label{(LPforTDIdef)}
\end{equation}
has an integral optimal solution $y$ for every integral vector 
\REV{
$c \in {\bf Z}\sp{S}$ 
for which the maximum is finite.  
More generally (see, \cite[Vol.~A, p.~77]{Schrijverbook}),
a rational linear system 
$[Q_{1} x \geq p_{1}, Q_{2} x =p_{2}]$ is defined to be {\bf TDI} 
if the system $[Q_{1} x\geq p_{1}, Q_{2} x \geq p_{2}, -Q_{2} x \geq -p_{2}]$ is TDI, 
which is equivalent to requiring that, for each integral vector 
$c \in {\bf Z}\sp{S}$, the dual of the primal linear program 
$\min\{cx:  Q_{1} x \geq p_{1}, Q_{2} x=p_{2} \}$ 
has an integer-valued optimum solution, if it has a finite optimum.
}

\REV{%
Edmonds and Giles \cite{Edmonds-Giles84} called a system
}%
 $Qx\geq p$ 
{\bf box-totally dual integral}  ({\bf box-TDI}) 
if the system
$[ Qx\geq p, f\leq x\leq g ]$ 
is TDI for every choice of rational
(finite-valued) bounding 
\REV{
vectors
}
$f\leq g$.  
This definition can be extended to linear systems including equations, as follows.  
\REV{
A linear system $[ Q_{1}x\geq p_{1}, Q_{2}x=p_{2} ]$ 
is called
 {\bf box-TDI} if the 
linear system
$[Q_{1} x\geq p_{1}, Q_{2} x \geq p_{2}, -Q_{2} x \geq -p_{2}]$ 
is box-TDI.
It follows from these definitions that
a linear system $[Q_1x\geq p_1, Q_2x=p_2]$ is box-TDI 
if and only if the system $[Q_1x\geq p_1, Q_2x =
p_2, f\leq x\leq g]$ is TDI for every choice of rational
(finite-valued) bounding functions $f\leq g$.  
}

A polyhedron is called a {\bf box-TDI polyhedron} if it can be
described by a box-TDI system.  
Edmonds and Giles proved basic properties of box-TDI systems,
while the paper of Cook \cite{Cook86} includes further important results on box-TDI polyhedra.  
For a rich overview of the topic, see the book of Schrijver \cite{Schrijverbook0}
and the recent paper of Chervet, Grappe, and Robert \cite{CGR}.
\REV{
The convex hull of four vectors
$(1,1,1,0,0,0)$,
$(1,0,0,1,0,0)$,
$(0,1,0,0,1,0)$, 
$(0,0,1,0,0,1)$
is a known example of a non-box-TDI (0,1)-polyhedron
(a face of the stable set polytope of a graph known as $S_{3}$) \cite{Cameron89,CGR}.
}

Our goal is to show that a result analogous to Theorem~\ref{min-max.orig} 
holds for the set $\odotZ{R}$ of integral elements of a box-TDI polyhedron $R$,
\REV{
that is, for a discrete box-TDI set.
}
An important special case is when $Q$ is a TU (totally unimodular) matrix.  
This includes the special case of L-convex or L$\sp{\natural}$-convex sets.  
It can be proved that L$_{2}\sp{\natural}$-convex (in particular, L$_{2}$-convex) 
sets are also discrete box-TDI sets.
Another special case is the one of integral submodular flows, 
in particular, M$_{2}$-convex and M$_{2}\sp{\natural}$-convex sets.

\subsection{Properties and operations}

In this section, we collect some basic properties of box-TDI systems
and polyhedra, which shall serve as useful tools for our later investigations.

\begin{proposition} [{\cite[Theorem 22.7]{Schrijverbook0}}] \label{alsotdi} 
A box-TDI system is TDI.  
\finbox 
\end{proposition}

\begin{proposition} [\cite{Cook86}] \label{TDI-box} 
Any TDI linear system defining a box-TDI polyhedron is box-TDI.  
\finbox 
\end{proposition}

Let $Qx\geq p$ be a box-TDI system and let $R:=\{x:  Qx\geq p\}$.  
For technical simplicity, we formulate the next propositions only 
for this form but emphasize that
each proposition below extends to the case when the system is given 
in the more general form $[ Q_{1}x\geq p_{1}, Q_{2}x=p_{2} ]$,
\REV{
which, by definition, is box-TDI if and only if 
$[Q_{1} x\geq p_{1}, Q_{2} x \geq p_{2}, -Q_{2} x \geq -p_{2}]$ is box-TDI.
}


\begin{proposition}[{\cite[Theorem 5.34]{Schrijverbook}}] \label{eltol} 
\REV{
For a rational vector $z\sp{*}$,  
}
let $p_{0}:=p-Qz\sp{*}$.  
Then the system $Qx\geq p_{0}$ is box-TDI.  
\end{proposition}
\Proof 
\RED{
(A proof is given here for completeness, as it is omitted in \cite{Schrijverbook}.)
Let $f$ and $g$ be any finite-valued rational bounding vectors 
with $f\leq g$.
Let $c$ be an integral vector for which the dual problem
\begin{equation}  \label{(transprf1)}
\ \max \{y p_{0} + uf - v g:   yQ + u -v =c, (y,u,v)\geq 0 \}
\end{equation}
has a finite optimal value. 
By using the definition $p_{0}=p-Qz\sp{*}$ and the constraint $yQ = c - u + v$,
we can rewrite the objective function in \eqref{(transprf1)} as
\begin{align*}
y p_{0} + uf - v g 
& = 
y (p-Qz\sp{*}) + uf - v g  
\nonumber \\ & =
y p- (c - u + v)z\sp{*} + uf - v g 
\nonumber \\ & =
y p + u(f+ z\sp{*}) - v (g+ z\sp{*}) - c z\sp{*} ,
\end{align*}
where the last term $c z\sp{*}$ is a constant independent of $(y,u,v)$.
Therefore, $(y,u,v)$ is an optimal solution to \eqref{(transprf1)} if and only if 
it is an optimal solution to 
\begin{equation}  \label{(transprf3)}
\ \max \{y p + u(f+ z\sp{*}) - v (g+ z\sp{*}):   yQ + u -v =c, (y,u,v)\geq 0 \} .
\end{equation}
Since the system 
$[ Qx\geq p, f+z\sp{*} \leq x\leq g+z\sp{*} ]$ is TDI
by the assumed box-TDI-ness of $Qx\geq p$,
the problem \eqref{(transprf3)} has 
an integral optimal solution $(y,u,v)$.
Therefore, the system $Qx\geq p_{0}$ is box-TDI.  
}
\finbox

\begin{proposition} \label{negating} 
If $Q'$ is a matrix obtained from $Q$ by negating some columns of $Q$, 
then the system $Q'x'\geq p$ is also box-TDI.  
\end{proposition}

\Proof 
It suffices to prove the special case when we negate the first column of $Q$.  
Let $Q'$ denote the matrix arising in this way.  
Let $f'\leq g'$ be rational bounding vectors and $c'$ an integer cost-function.  
We have to show that the dual program
\begin{equation}  
\max \{ yp + f'u - g'v:  yQ'+u-v = c', (y,u,v)\geq 0\}
\label{(nega)} 
\end{equation} 
has an integral optimal solution $(y,u,v)\geq 0$.
Let $c$ denote the vector obtained from $c'$ by negating its first component.  
Let $f$ denote the vector obtained from $f'$ 
by replacing its first component $f'(1)$ to $-g'(1)$, 
and let $g$ denote the vector obtained from $g'$ 
by replacing its first component $g'(1)$ to $-f'(1)$.  
Then $yQ'+u'-v' = c'$ if and only if $yQ+u-v = c$,
where $(u',v')$ arises from $(u,v)$ by interchanging their first components.
Furthermore, 
$yp + f'u' - g'v'= yp +fu-gv$.  
By the box-TDI-ness of the system $Qx\geq p$, 
there is an integer-valued optimal solution $(y,u,v)$ to
\begin{equation}  
\max \{ yp + fu - gv:  yQ+u-v = c, (y,u,v)\geq 0\} 
\label{(nega2)}
\end{equation}
\noindent 
and hence $(y,u',v')$ is an integer-valued optimal solution to \eqref{(nega)}.  
\finbox

\begin{proposition} [{\cite[p.~323]{Schrijverbook0}}] \label{deletion} 
The system obtained from a box-TDI system $Qx\geq p$ by deleting some columns of $Q$ is box-TDI.  
\finbox 
\end{proposition}

\begin{proposition} [\cite{Edmonds-Giles84}] \label{dupli} 
If $Q'$ is a matrix obtained from $Q$ by duplicating some columns of $Q$, 
then the system $Q'x'\geq p$ is also box-TDI.  
\finbox 
\end{proposition}

\begin{proposition} [{\cite[p.~323]{Schrijverbook0}}] \label{projection} 
The projection of a box-TDI polyhedron along a coordinate axis is box-TDI.
\finbox 
\end{proposition}

\begin{proposition} \label{superfluous} 
Let $Qx\geq p$ be a (box-) TDI system defining
\REV{
the
}
 polyhedron $R:=\{x:Qx\geq p\}$.  
Let $qx\geq \beta $ be an inequality which is superfluous 
in the sense that every member $x$ of $R$ satisfies $qx\geq \beta $. 
Then the system $[ Qx\geq p, qx\geq \beta ]$ is also (box-) TDI.  
\end{proposition}

\Proof 
Let $c$ be an integral cost-function and let $y_{0}$ 
be an integral dual optimum ensured by the TDI-ness of $Qx\geq p$.  
Since by adding a superfluous inequality to a linear system 
does not change the primal optimum value, 
the dual optimum value does not change either.
Therefore, by extending $y_{0}$ by a new zero-valued dual component
corresponding to the primal inequality $qx\geq \beta $, 
we obtain an integral dual solution $(y_{0},0)$ to the dual of the primal problem
$\min \{cx:  Qx\geq p, \  qx\geq \beta \}$.

The statement for box-TDI-ness follows from the first part since if
$qx\leq \beta $ is superfluous with respect to the system $Qx\geq p$,
then it is superfluous, as well, for the system 
$[ Qx\geq p, \  f\leq x\leq g ]$ 
for any pair of bounding functions $f\leq g$.  
\finbox

\begin{proposition} \label{unbo} 
Let $Qx\geq p$ be a box-TDI system.  
Let
$f':S\rightarrow {\bf Q} \cup \{-\infty \}$ 
and $g':S\rightarrow {\bf Q} \cup \{+\infty \}$ 
be rational bounding vectors with $f'\leq g'$.
Then $[ Qx\geq p,  \  f'\leq x\leq g' ]$ 
is also box-TDI, and (hence) TDI.
\end{proposition}

\Proof 
We have to show for any choice $f:S\rightarrow {\bf Q}$ 
and
$g:S\rightarrow {\bf Q}$
of finite-valued rational bounds that the system 
\begin{equation}  
  [ Qx\geq p, \  f'\leq x\leq g', \  f\leq x\leq g ]
 \label{(fgfg)} 
\end{equation} 
is TDI.  
Let $f_{0}$ be the componentwise maximum of $f$ and $f'$, 
and let $g_{0}$ be the componentwise minimum of $g$ and $g'$.  
Then $f_{0}$ and $g_{0}$ are finite-valued and hence the system 
\begin{equation}   
[ Qx\geq p,  \ f_{0}\leq x\leq g_{0} ]
\label{(f0g0)} 
\end{equation} 
\RED{
is TDI since $Qx\geq p$ is box-TDI.
}
Since the system in \eqref{(fgfg)} arises from the system in \eqref{(f0g0)} by
adding superfluous inequalities, Proposition~\ref{superfluous} implies
that \eqref{(fgfg)} is indeed TDI, as required. 
\finbox.

\medskip

\begin{proposition} \label{ezkell} 
Let $Qx\geq p$ be a box-TDI system defining the box-TDI polyhedron $R:=\{x:  Qx\geq p\}$, 
let  $z\sp{*}\in \odotZ{R}$, and $p_{0}:=p-Qz\sp{*}$.  
Then the system 
\begin{equation}
[ (x_{1},x_{2}) \geq 0, \ Qx_{2} - Qx_{1} \geq p_{0}  ] 
\label{(boxtdi1)} 
\end{equation}
is box-TDI.  
\end{proposition}

\Proof 
By Proposition~\ref{eltol}, the system $Qx\geq p_{0}$ is box-TDI.
By Proposition~\ref{dupli}, $Qx_{2} + Qx_{1} \geq p_{0}$ is box-TDI.  
By applying Proposition~\ref{negating} to the matrix $(Q,Q)$, we get that
$Qx_{2} - Qx_{1} \geq p_{0}$ is box-TDI.  
And finally, by Proposition~\ref{unbo}, 
the system $[ (x_{1},x_{2}) \geq 0, \ Qx_{2} - Qx_{1} \geq p_{0} ]$ is box-TDI.  
\finbox


\medskip

\REV{
 A polyhedron is called {\bf box-integer} \cite{CGR,Schrijverbook} 
if its intersection with any integral box is integral.  
For a positive integer $k$ the {\bf $k$-dilation} \  $kR$ \  of a polyhedron 
$R=\{x:Qx\geq p\}$ is defined by $\{x:  Qx\geq kp\}$.  
Any $k$-dilation is called an (integer) {\bf dilation} of $R$.
}

\REV{
\begin{proposition} [\cite{CGR}] \label{dilat} 
An integer polyhedron $R$ is box-TDI 
if and only if each of its integer dilation is box-integer.
\finbox 
\end{proposition}
}

\REV{
\begin{remark} \rm \label{RMcautionTDI}
In this section, we have indicated that
some natural basic operations preserve \hbox{(box-)} total dual integrality.
It should, however, be remarked that one has to be cautious in
formulating such results since there are other ``natural'' operations
that do not preserve (box-) TDI-ness.  
For example, a remark of Schrijver's book \cite[p.~323]{Schrijverbook0} 
cites a counter-example of Cook \cite{Cook83} which demonstrates that the statement in
Proposition \ref{ezkell} does not hold anymore if we replace box-TDI-ness by TDI-ness.  
Another negative result is that the
TDI-ness of the system $[Qx\geq p_{1}, Qx\geq p_{2}]$ does not imply the
TDI-ness of the system $Qx\geq p_{1}+p_{2}$.
Also, R. Grappe pointed out
that adding a unit vector $(1,0,0,\dots ,0)$ as a column to the
constraint matrix in a box-TDI system may destroy box-TDI-ness.
\hfill $\bullet$ 
\end{remark}  
}

\subsection{The main tool}

The following result is the main tool in proving the min-max theorem
in Section \ref{main-min-max}.

\begin{theorem} \label{inverse} 
Let $Q$ be an integral matrix, $p$ an integral vector, 
and suppose that the linear system $Qx\geq p$ is box-TDI.  
Let $z\sp{*}$ be an integral element of 
\REV{
the
}
polyhedron 
$R:=\{x: Qx\geq p\} \subseteq {\bf R}\sp{S}$, 
and let $\ell:  S\rightarrow {\bf Z} \cup \{-\infty \}$
and $u:  S\rightarrow {\bf Z} \cup \{+\infty \}$
be integer-valued bounding vectors on $S$ for which $\ell\leq u$.
There exists an integer-valued non-negative vector $y\sp{*}$ such that
$\ell\leq y\sp{*}Q \leq u$ and $y\sp{*}(Qz\sp{*}-p) = 0$ 
if and only if
\begin{equation} 
 \widetilde \ell(S\sp{-}) \leq \widetilde u(S\sp{+}) 
 \label{(fSgS)}
\end{equation}
\noindent 
holds for every pair 
\RED{
$( S\sp{-},S\sp{+} )$ 
}
of disjoint subsets of $S$ for which
\begin{equation} 
 z':= z\sp{*} + \chi_{S\sp{+}} - \chi_{S\sp{-}} \in R,
\label{(mprime)} 
\end{equation} 
\REV{
where $\chi_{S\sp{+}}$ and $\chi_{S\sp{-}}$ denote
the characteristic vectors of $S\sp{+}$ and $S\sp{-}$, respectively.
}
\end{theorem}

\Proof 
Necessity of \eqref{(fSgS)}.  
Let $y\sp{*}$ be a function meeting the requirements, 
$w\sp{*}:=y\sp{*}Q$, and let
\RED{
$( S\sp{-},S\sp{+} )$ 
}
be a pair meeting \eqref{(mprime)}.  
\REV{
Then, by complementary slackness of the pair of linear programs \eqref{(LPforTDIdef)}
for $c = w\sp{*}$,
}
$y\sp{*}(Qz\sp{*}-p) = 0$ implies that 
$z\sp{*}$ is $w\sp{*}$-minimizer of $R$, and hence 
\[
w\sp{*}z\sp{*} \leq w\sp{*}z' 
= w\sp{*}z\sp{*} + \widetilde w\sp{*}(S\sp{+}) - \widetilde w\sp{*}(S\sp{-}) 
\leq w\sp{*}z\sp{*} + \widetilde u(S\sp{+}) - \widetilde \ell(S\sp{-}),
\] 
from which \eqref{(fSgS)} follows.  (Note that
$\widetilde u(S\sp{+})=+\infty $ and $\widetilde \ell(S\sp{-})=-\infty $
may occur.)

Sufficiency of \eqref{(fSgS)}.  
Let $p_{0}:=p-Qz\sp{*}$.  By the linear programming duality theorem, we have
\begin{align}  
& \min \{ux_{2} - \ell x_{1}:  (x_{1},x_{2})\geq 0, \ Qx_{2} - Qx_{1} \geq p_{0}\}
\label{(primal1)} 
\\
& = \max \{yp_{0} :  y\geq 0, \ yQ \leq u, \ y(-Q) \leq -\ell\}.
\label{(dual1)} 
\end{align}
\noindent 
Formally, this is correct only if both $u$ and $\ell$ are
finite-valued.  To get the right pair of dual programs for the general
case, one must remove the columns of $Q$ corresponding to elements $s$
with $u(s)=+\infty $ and remove the columns of $-Q$ corresponding to
elements $s$ with $\ell(s)=-\infty $. But in order to avoid notational
difficulties, with this remark in mind, we work with the dual linear
programs \eqref{(primal1)} and \eqref{(dual1)}.

By Proposition~\ref{ezkell}, the linear system in \eqref{(primal1)} is box-TDI.  
Let $M$ denote the common optimum value of the primal and
the dual programs.  Since $y\geq 0$ and $p_{0}\leq 0$, we have $M\leq 0$.

\begin{claim}  \label{Mnonneg} 
 $M=0$.  
\end{claim}

\Proof 
Suppose indirectly that $M<0$.  
Then there is a solution $(x'_{1},x'_{2})$ to \eqref{(primal1)} 
for which $ux'_{2} - \ell x'_{1} = M < 0$.
By the definition of $p_{0}$, the primal constraint 
$Qx'_{2} - Qx'_{1} \geq p_{0}$ is equivalent to 
\REV{
$z_{1}\sp{*}:=z\sp{*} + x'_{2} -x'_{1}\in R$.  
}
Since both $z\sp{*}$ and 
\REV{
$z_{1}\sp{*}$  
}
are in $R$, 
the line segment connecting $z\sp{*}$ and 
\REV{
$z_{1}\sp{*}$  
}
also lies in $R$, 
that is, for any $\varepsilon $ with $0\leq \varepsilon \leq 1$, 
the vector $z\sp{*} + \varepsilon (x'_{2} - x'_{1})$ belongs to $R$, 
or equivalently $\varepsilon (Qx'_{2}-Qx'_{1})\geq p_{0}.$
We can choose $\varepsilon $ in such a way that $0<\varepsilon \leq 1$,
\begin{equation}  
x''_{1}(s):= \varepsilon x'_{1}(s) \leq 1 \ \hbox{ and }\ \
x''_{2}(s):=\varepsilon x'_{2}(s) \leq 1 \ \hbox{ for every }\  s\in S.
\label{(01box)} 
\end{equation}
\noindent 
Clearly, $Qx''_{2}-Qx''_{1}\geq p_{0}$ and
\begin{equation}  
ux''_{2}-\ell x''_{1} = \varepsilon (ux'_{2}-\ell x'_{1}) = \varepsilon M <0.  
\label{(negat)} 
\end{equation}
\noindent 
These imply that the linear system 
$[ (x_{1},x_{2})\geq 0, \  Qx_{2}-Qx_{1}\geq p_{0} ]$ 
in \eqref{(primal1)} has a solution meeting \eqref{(01box)} and \eqref{(negat)}.  
The box total dual integrality of the linear system in \eqref{(primal1)} 
implies that there is
\REV{
a
}
$\{0,1\}$-valued solution $(x\sp{*}_{1},x\sp{*}_{2})$ 
for which $M\sp{*}:=ux\sp{*}_{2}-\ell x\sp{*}_{1} <0$.

Furthermore, we can also assume that no element $s\in S$ exists with
$x\sp{*}_{1}(s)=1= x\sp{*}_{2}(s)$ 
since in this case we could reduce both values by 1, 
and then $\ell(s)\leq u(s)$ would imply for the revised
$(x\sp{*}_{1},x\sp{*}_{2})$ 
that $ ux\sp{*}_{2} - \ell x\sp{*}_{1} = M\sp{*}- u(s) + \ell(s) \leq M\sp{*} <0$.

Let $S\sp{+}:=\{s\in S:  x\sp{*}_{2}(s)>0\}$ 
and $S\sp{-}:=\{s\in S:  x\sp{*}_{1}(s)>0\}$.  
Then $S\sp{+}$ and $S\sp{-}$ are disjoint for which 
$\widetilde u(S\sp{+})= ux\sp{*}_{2} < \ell x\sp{*}_{1} = \widetilde \ell (S\sp{-})$, 
contradicting \eqref{(fSgS)}. 
\finbox

\medskip

As $M=0$, the box-TDI-ness of the linear system in \eqref{(primal1)}
implies that the dual problem in \eqref{(dual1)} has an integer-valued
solution $y\sp{*}$ for which $y\sp{*}p_{0}=M=0$, that is, 
$y\sp{*}(Qz\sp{*}-p)=0$.  
Furthermore $\ell\leq w\sp{*} \leq u$ holds for $w\sp{*}:=y\sp{*}Q$, as required.  
\BB

\begin{corollary} \label{inverz2} 
Let $Q, p, R, \ell, u$, and $z\sp{*}$ be the
same as in Theorem {\rm \ref{inverse}}.  
There exists an integer-valued cost-function $w\sp{*}$ on $S$ 
for which $\ell\leq w\sp{*} \leq u$ and
$z\sp{*}$ is a $w\sp{*}$-minimizer of $R$ 
if and only if 
\eqref{(fSgS)} holds for every pair 
\RED{
$( S\sp{-},S\sp{+} )$ 
}
of disjoint subsets of $S$ meeting \eqref{(mprime)}.  
\end{corollary}

\Proof 
The corollary follows immediately from Theorem~\ref{inverse}
once we make the standard observation from linear programming 
that a primal solution $z\sp{*}$ is a $w\sp{*}$-minimizer of $R$ 
if and only if
there is a dual solution $y\sp{*}$ meeting the optimality criteria,
that is, 
$y\sp{*} \geq 0$, $y\sp{*}Q=w\sp{*}$, and $y\sp{*}(Qz\sp{*}-p) = 0$.  
\finbox


\section{Min-max theorem for $\Phi$} 
\label{main-min-max}

\subsection{Preparation}

Let $\varphi :{\bf Z} \rightarrow {\bf Z} \cup \{+\infty \}$ 
be an arbitrary integer-valued function on ${\bf Z}$ 
allowing the $+\infty$ value.  
We say that an ordered pair 
\REV{
$( k\sp{*},\ell\sp{*} )$ 
}
of integers is 
{\bf $\varphi$-fitting} if 
\begin{equation}  
\varphi (k\sp{*}) - \varphi(k\sp{*}-1)  \leq \ell\sp{*} 
\leq \varphi (k\sp{*}+1) - \varphi (k\sp{*})
\label{(fitting1)} 
\end{equation}
\noindent 
or more concisely
\begin{equation}  
\varphi '(k\sp{*}-1) \leq \ell\sp{*} \leq \varphi '(k\sp{*}).
\label{(fitting2)} 
\end{equation}

Let $\Phi $ be a separable function on ${\bf Z}\sp{S}$ defined by
univariate integer-valued functions $\varphi_{s}$ ($s\in S$).  
We say that an ordered pair 
\REV{
$( z\sp{*},w\sp{*} )$ 
}
of vectors from ${\bf Z}\sp{S}$ 
is {\bf $\Phi$-fitting} 
if 
\REV{
$( z\sp{*}(s),w\sp{*}(s) )$ 
}
is $\varphi_{s}$-fitting for each $s\in S$, that is,
\begin{equation} 
 \varphi_{s}(z\sp{*}(s)) - \varphi_{s}(z\sp{*}(s)-1) \leq w\sp{*}(s)
\leq \varphi_{s}(z\sp{*}(s)+1) - \varphi_{s}(z\sp{*}(s))\ \ \hbox{for every}\ \ s\in S, 
\label{(fitting)} 
\end{equation}
\noindent 
which can concisely be written as follows:
\begin{equation}  
\Phi' (z\sp{*}- {\bf \underline{1}} ) \ \leq \ w\sp{*} \ \leq \ \Phi' (z\sp{*}).  
\label{(fittingxx)} 
\end{equation}

As a preparation, we need the following proposition.

\begin{proposition} \label{fitting} 
Let $\varphi$ be an integer-valued discrete convex function 
and let 
\REV{
$( k\sp{*},\ell\sp{*} )$ 
}
be a $\varphi$-fitting pair of integers.  
Then 
\begin{equation} 
 \ell\sp{*}k\sp{*} - \varphi (k\sp{*})  \geq \ell\sp{*}k - \varphi (k) 
\quad \hbox{\rm for every integer}\ \ k
\label{(fitting3)} 
\end{equation} 
\noindent 
(or equivalently $\varphi
\sp{\bullet}(\ell\sp{*})= \ell\sp{*}k\sp{*} - \varphi (k\sp{*})$ where $\varphi
\sp{\bullet}$ denotes the discrete conjugate of $\varphi$).  
\end{proposition}

\Proof 
Suppose indirectly that there is an integer $k_{0}$ for which 
\begin{equation}  
 \ell\sp{*}k\sp{*} - \varphi (k\sp{*}) < \ell\sp{*}k_{0} - \varphi (k_{0}).
\label{(indir)} 
\end{equation}
\noindent 
If $k_{0}>k\sp{*}$, we may assume that $k_{0}$ is minimal, and hence 
\begin{equation}  
\ell\sp{*}k\sp{*} - \varphi (k\sp{*}) \geq \ell\sp{*}(k_{0}-1) - \varphi (k_{0}-1).  
\label{(k0min)} 
\end{equation} 
By subtracting \eqref{(k0min)} from \eqref{(indir)}, we get 
$0< \ell\sp{*} - (\varphi (k_{0}) - \varphi (k_{0}-1))$.  
This and the convexity of $\varphi$ imply that 
$\ell\sp{*} > \varphi (k_{0}) - \varphi (k_{0}-1)\geq \varphi (k\sp{*}+1) - \varphi (k\sp{*})$, 
in contradiction to the second inequality in \eqref{(fitting1)}.

Analogously, if $k_{0}<k\sp{*}$, we may assume that $k_{0}$ is maximal, and hence 
\begin{equation}  
\ell\sp{*}k\sp{*} - \varphi (k\sp{*}) \geq \ell\sp{*}(k_{0}+1) -
\varphi (k_{0}+1).  
\label{(k0max)} 
\end{equation} 
By subtracting \eqref{(indir)} from \eqref{(k0max)}, we get 
$0 > \ell\sp{*} - (\varphi (k_{0}+1) - \varphi (k_{0}))$.  
This and the convexity of $\varphi$ imply that
$\ell\sp{*} < \varphi (k_{0}+1) - \varphi (k_{0})\leq \varphi (k\sp{*}) - \varphi (k\sp{*}-1)$, 
in contradiction to the first inequality in \eqref{(fitting1)}.  
\finbox

\begin{remark} \rm
There is a standard concept and terminology in (discrete) convex analysis 
that is equivalent in the present case to $\varphi$-fitting.  
Namely, 
\REV{
$\ell \sp{*}$ satisfying the condition \eqref{(fitting2)} is
called a subgradient 
of $\varphi$ at $k\sp{*}$,
and the set of 
these subgradients 
is called the subdifferential of $\varphi$ at $k\sp{*}$,
usually denoted by $\partial \varphi (k\sp{*})$.
}
Therefore, 
\REV{
$( k\sp{*},\ell\sp{*} )$ 
}
is $\varphi$-fitting 
if and only if 
$\ell \sp{*}\in \partial \varphi (k\sp{*})$.  
Proposition~\ref{fitting} is a restatement of the well-known fact 
that 
$\varphi(k\sp{*}) + \varphi \sp{\bullet}(\ell\sp{*}) = k\sp{*} \ell\sp{*}$ 
holds if and only if 
$\ell \sp{*}\in \partial \varphi (k\sp{*})$.  
\hfill $\bullet$ 
\end{remark}

\subsection{Main results}

Let 
\REV{
$R =\{x:  Qx\geq p\} \subseteq {\bf R}\sp{S}$ 
}
be an arbitrary integral polyhedron and
$z\sp{*}$ an element of $\odotZ{R}$.  
Let $\varphi_{s}$ be an integer-valued discrete convex function on ${\bf Z}$ 
for each $s\in S$ and let $\Phi$ denote the separable discrete convex function
defined in \eqref{(Phidef)} by the univariate functions $\varphi_{s}$ $(s\in S)$.

Let $y\sp{*}$ be a vector whose components correspond to the rows of $Q$.  
We say that the ordered pair 
\RED{
$( z\sp{*},y\sp{*} )$ 
}
of integral vectors is 
{\bf $\Phi$-compatible} with respect to $Q$ 
(or, shortly $\Phi$-compatible) 
if 
\REV{
$( z\sp{*},w\sp{*} )$ 
}
 is $\Phi$-fitting where $w\sp{*}:=y\sp{*}Q$, that is,
\begin{equation}  
 \Phi' (z\sp{*}- {\bf \underline{1}} ) \ \leq \ y\sp{*}Q \ \leq \ \Phi' (z\sp{*}).  
\label{(fittingx)} 
\end{equation}

\begin{remark} \rm \label{Philinear.1} 
In the special case when $\Phi$ is a linear function, 
that is, $\Phi(z)=cz$ for a given vector $c\in {\bf Z}\sp{S}$, 
one has $\Phi'(z)=c$ for every $z\in {\bf Z}\sp{S}$.
Therefore, in this case, $\Phi$-compatibility given in \eqref{(fittingx)} 
is equivalent to $c\leq y\sp{*}Q\leq c$, that is,
$c=y\sp{*}Q$.  
\hfill $\bullet$ 
\end{remark}  

\begin{lemma} \label{lowerb} 
Let $\Phi$ be an integer-valued separable discrete convex function on ${\bf Z} \sp{S}$.  
Suppose for  $z\sp{*}\in \odotZ{R}$ and $y\sp{*} \geq 0$ 
that the pair 
\RED{
$( z\sp{*},y\sp{*} )$ 
}
is $\Phi$-compatible.  Then
\begin{equation}
  \Phi(z) \geq \Phi(z\sp{*}) - y\sp{*}(Q z\sp{*}- p) 
\label{(lowerb)}
\end{equation}
holds for every $z\in \odotZ{R}$.  
\end{lemma}

\Proof Let $w\sp{*}:=y\sp{*}Q$.  Since 
\REV{
$( z\sp{*},y\sp{*} )$ 
}
is $\Phi$-compatible, 
\REV{
$( z\sp{*},w\sp{*} )$ 
}
is a $\Phi$-fitting pair, and
we can apply Proposition~\ref{fitting} to $\varphi :=\varphi_{s}$,
$k\sp{*}:=z\sp{*}(s)$, $\ell\sp{*}:  = w\sp{*}(s)$, and $k:=z(s)$:
\begin{equation*}  
w\sp{*}(s)z\sp{*}(s) - \varphi_{s}(z\sp{*}(s)) \geq w\sp{*}(s)z(s) - \varphi_{s}(z(s)), 
\end{equation*} 
that is,
\begin{equation} 
 \varphi_{s}(z(s)) \geq 
  w\sp{*}(s)z(s) - [ w\sp{*}(s)z\sp{*}(s) - \varphi_{s}(z\sp{*}(s))].  
\label{(zpi2)} 
\end{equation}
\REV{
On the other hand, we have
}
\begin{equation} 
\sum_{s\in S}w\sp{*}(s)z(s) = w\sp{*} z = (y\sp{*}Q)z = y\sp{*}(Qz) \geq y\sp{*}p
 \label{(zpi3)} 
\end{equation}
\REV{
since $z \in \odotZ{R}$ and $y\sp{*} \geq 0$.
It follows from \eqref{(zpi2)} and \eqref{(zpi3)} that   
}
\begin{align} 
 \Phi(z) & = \sum_{s\in S} \varphi_{s}(z(s))
\nonumber \\ & 
 \geq \sum_{s\in S}w\sp{*}(s)z(s) 
 - \bigg[ \sum_{s\in S}w\sp{*}(s)z\sp{*}(s) - \sum_{s\in S}\varphi_{s}(z\sp{*}(s)) \bigg] 
\nonumber \\ &
\geq 
 y\sp{*}p - \bigg[ \sum_{s\in S}w\sp{*}(s)z\sp{*}(s) - \sum_{s\in S}\varphi_{s}(z\sp{*}(s)) \bigg] 
\nonumber \\ & 
=  y\sp{*}p - [(y\sp{*}Q) z\sp{*} - \Phi (z\sp{*})] 
 =
\Phi(z\sp{*}) - y\sp{*}(Q z\sp{*}- p), 
 \label{(zpi4)} 
\end{align} 
as required. 
\finbox

\medskip

The new min-max theorem for the case when $R$ is an integral box-TDI polyhedron is as follows.

\begin{theorem} \label{min-max-uj} 
Let $\varphi_{s}:{\bf Z} \rightarrow {\bf Z} \cup \{+\infty \}$ 
be an integer-valued discrete convex function on ${\bf Z}$ for each $s\in S$ and 
let $\Phi$ denote the separable discrete convex function 
defined by the univariate functions $\varphi_{s}$ $(s\in S)$.  
Suppose for an integral matrix $Q$ and 
\REV{
an
}
integral vector $p$ 
that $Qx\geq p$ is a box-TDI system defining a non-empty integral (box-TDI) polyhedron 
$R:=\{x:  Qx\geq p\}\subseteq {\bf R}\sp{S}$ 
such that $\Phi$ is finite-valued on $\odotZ{R}$.  
Then $\Phi$ is bounded from below on $\odotZ{R}$ 
if and only if 
there exists an element $z\in \odotZ{R}$ and an integral vector $y\geq 0$ for which
\REV{
$( z,y )$ 
}
is $\Phi$-compatible with respect to $Q$.  
Moreover, if $\Phi$ is bounded from below on $\odotZ{R}$, 
then the following min-max formula holds:
\begin{align}  
 & \min \{ \Phi (z) :  z\in \odotZ{R} \}  
\nonumber \\ & =
  \max \{ \Phi(z) - y(Qz-p) :  \ z\in \odotZ{R} , \ y\geq 0 \ \
  \hbox{\rm integer-valued, \ 
\REV{
$( z,y )$ 
}
\ $\Phi$-compatible} \} .
\label{(dualmaxuj)} 
\end{align}
\noindent
In addition, an optimal vector $y\sp{*}$ in \eqref{(dualmaxuj)} 
can be chosen in such a way that the number of
its positive components is at most 
$2\vert S\vert$.  
\end{theorem}

\Proof 
Suppose first that there is a 
$\Phi$-compatible pair 
\REV{
$( z\sp{*},y\sp{*} )$ 
}
\REV{
with $z\sp{*} \in \odotZ{R}$ and $y\sp{*} \geq 0$. 
}
Then Lemma \ref{lowerb} implies that $\Phi$ is bounded from below and that 
$\min \geq \max$.

Suppose now that $\Phi$ is bounded from below on $\odotZ{R}$.  Since
$\Phi$ is integer-valued, $\odotZ{R}$ has a $\Phi$-minimizer element
$z\sp{*}$.  We are going to show that there is an integer-valued vector
$y\sp{*} \geq 0$ for which the following optimality criteria hold:
\begin{align}  
& y\sp{*}(Qz\sp{*}-p) = 0, 
\label{(optcrit1bD)} 
\\ &
 \Phi' (z\sp{*}- {\bf \underline{1}} ) \leq y\sp{*}Q \leq \Phi' (z\sp{*}).
\label{(optcrit2bD)} 
\end{align}
\noindent 
This will imply that
a $\Phi$-compatible pair in question indeed exists 
\REV{
which shows the equality in \eqref{(dualmaxuj)}.  
}


Define bounding vectors $\ell$ and $u$ on $S$, as follows.  
For $s\in S$, let
\[ 
\ell(s):=\varphi_{s}'(z\sp{*}(s)-1) \quad \hbox{ and } \quad 
u(s):= \varphi_{s}'(z\sp{*}(s)),
\]
where $\ell(s)$ may be $-\infty $ and $u(s)$ may be $+\infty$. 
The discrete convexity of $\varphi_{s}$ implies that $\ell(s)\leq u(s)$.  
Note that 
\REV{
$( z\sp{*},y\sp{*} )$ 
}
is $\Phi$-compatible precisely if 
$\ell \leq y\sp{*}Q \leq u$.

\begin{claim} 
The inequality 
$\widetilde \ell(S\sp{-}) \leq \widetilde u(S\sp{+})$ 
in \eqref{(fSgS)} holds for every pair 
\RED{
$(S\sp{-},S\sp{+})$ 
}
of disjoint subsets 
of $S$ for which 
$z':= z\sp{*} + \chi_{S\sp{+}} - \chi_{S\sp{-}} \in R$. 
\end{claim}

\Proof 
As $z\sp{*}$ is a $\Phi$-minimizer, we have 
$\Phi(z\sp{*}) \leq \Phi(z')$.  
Furthermore
\begin{align*}
\Phi(z') &= \sum_{s\in S}\varphi_{s}(z'(s))
\\ &
 = \sum_{s\in S-(S\sp{+} \cup S\sp{-})}\varphi_{s}(z\sp{*}(s)) 
  + \sum_{s\in S\sp{+}}\varphi_{s}(z\sp{*}(s)+1) 
  + \sum_{s\in S\sp{-}}\varphi_{s}(z\sp{*}(s)-1) 
\\ &
= \sum_{s\in S}\varphi_{s}(z\sp{*}(s)) 
  + \sum_{s\in S\sp{+}} \big[\varphi_{s}(z\sp{*}(s)+1) -\varphi_{s}(z\sp{*}(s)) \big] 
  - \sum_{s\in S\sp{-}} \big[\varphi_{s}(z\sp{*}(s))- \varphi_{s}(z\sp{*}(s)-1) \big] 
\\ &
= \Phi(z\sp{*}) + \sum_{s\in S\sp{+}}\varphi '_{s}(z\sp{*}(s)) 
- \sum_{s\in S\sp{-}}\varphi '_{s}(z\sp{*}(s)-1) 
\\ &
=
\Phi(z\sp{*}) + \widetilde u(S\sp{+}) - \widetilde \ell(S\sp{-}) 
\\ &
\leq
\Phi(z') + \widetilde u(S\sp{+}) - \widetilde \ell(S\sp{-}), 
\end{align*}
from which $\widetilde \ell(S\sp{-}) \leq \widetilde u(S\sp{+})$, as required.  
\finbox

\medskip

Theorem~\ref{inverse} implies the existence of the requested $y\sp{*}$
\REV{
satisfying \eqref{(optcrit1bD)} and \eqref{(optcrit2bD)}.
Since a box-TDI linear system is totally dual integral by 
Proposition~\ref{alsotdi},
the last statement about the number of positive components is 
a consequence of a theorem of
Cook, Fonlupt, and Schrijver \cite{CFS}
(see also Theorem 5.30 in the book of Schrijver \cite{Schrijverbook}).
}
\BB

\noindent 
\begin{remark} \rm  \label{RMhilbbas}
Cook, Fonlupt, and Schrijver \cite{CFS}
actually proved a slightly better bound $2\vert S\vert -1$ for the
number of non-zero variables, and this was later improved to 
$2\vert S\vert -2$ by Seb{\H o} \cite{Sebo90}.  
The point in Theorem \ref{min-max-uj}
(and its consequences below) is that there is a reasonably small bound.
\hfill $\bullet$ 
\end{remark}

\begin{remark} \rm 
Note that the dual objective function in \eqref{(dualmaxuj)} can be rewritten, as follows:
\begin{equation}
  \Phi(z) - y(Qz-p) = yp - [ (yQ)z - \Phi(z)].  
\label{(rewrit)}
\end{equation}
\noindent 
Furthermore, the characterization of boundedness in Theorem~\ref{min-max-uj} 
can be interpreted as a special case of the min-max formula 
when the minimum 
\REV{
in \eqref{(dualmaxuj)}
}
is $-\infty$ 
and the maximum in \eqref{(dualmaxuj)}, 
when taken over the empty set, is defined to be $-\infty $. 
Therefore, in the variations and applications of Theorem~\ref{min-max-uj} below, 
we shall not explicitly formulate the condition for the lower boundedness of $\Phi$.  
\hfill $\bullet$ 
\end{remark}

\begin{remark} \rm 
At first sight, this min-max theorem looks a bit strange
in the sense that in the maximization part, not only the usual dual
variable $y$ appears but integral members $z$ of the primal polyhedron
$R$ also show up.  Still, this form may be viewed as a proper min-max
theorem since the right-hand side is a straightforward lower bound for
the minimum, and for given $z\sp{*}$ and $y\sp{*}$, the validity of
optimality criteria \eqref{(optcrit1bD)} and \eqref{(optcrit2bD)} is easily checkable.  
\REV{
It is also noted that, apart from integrality, the min-max formula in \eqref{(dualmaxuj)} 
can be viewed as a variant of the Lagrangian duality as follows.
Let 
\[
L(x,y) :=
\begin{cases} 
\Phi(x) - y (Q x -p) & \quad \hbox{\rm if}\quad  y \geq 0, 
\cr 
-\infty & \quad \hbox{\rm otherwise} ,
\end{cases}
\]
which is the standard Lagrangian function
for the minimization of $\Phi(x)$ subject to the constraint $Q x -p \geq 0$.  
Then the Lagrangian dual problem is to maximize $\Psi(y) = \min_{x} L(x,y)$ over all $y \geq 0$.
When $\Phi(x)$ is convex, the minimum of $L(x,y)$ with respect to $x$ is attained by $x$ at which
$y Q$ is a subgradient of $\Phi$, that is, $yQ \in \partial \Phi(x)$.
Thus the dual problem
$\max \{ \Psi(y) : \ y\geq 0 \}$
 may be written as
$\max \{ \Phi(x) - y(Qx -p) :  \ y\geq 0, \  yQ \in \partial \Phi(x) \}$.
The constraint $yQ \in \partial \Phi(x)$ here is equivalent to
saying, in our present terminology, that 
\REV{
$( x,y )$ 
}
is $\Phi$-compatible.
Our dual problem in \eqref{(dualmaxuj)} is obtained by adding the constraint $x \in R$
to this Lagrangian dual problem.
}
\hfill $\bullet$ 
\end{remark}

\begin{remark} \rm \label{Philinear.2} 
In the special case when $\Phi(z)=cz$, the compatibility of $z$ and $y$, 
as observed in Remark~\ref{Philinear.1}, is equivalent to $yQ=c$.  
Furthermore, the dual objective function in \eqref{(dualmaxuj)} is as follows:
\[ 
\Phi(z) - y(Qz-p) = yp - [ (yQ)z - \Phi(z) ] = yp - [cz-cz] = yp, 
\]
showing that in this case we are back at the integral version of the
linear programming duality theorem formulated for box-TDI polyhedra.
\hfill $\bullet$ 
\end{remark}

It is useful to formulate separately the optimality criteria appearing
in \eqref{(optcrit1bD)} and \eqref{(optcrit2bD)}.

\begin{corollary} [Optimality criteria] \label{Optcrit} 
An element $z\sp{*}\in \odotZ{R}$ is a $\Phi$-minimizer 
if and only if 
there exists a non-negative integer-valued vector $y\sp{*}$ 
meeting the optimality criteria in \eqref{(optcrit1bD)} and \eqref{(optcrit2bD)}.  
\finbox
\end{corollary}

\subsection{Using discrete conjugate}

The min-max formula for the minimum of $\Phi$ can be described in a
more concise way in term of discrete conjugates.  
To this end, we need some easy observations.  
In Proposition~\ref{fitting}, we proved for a univariate discrete convex function $\varphi$ 
that if 
\REV{
$( k\sp{*},\ell\sp{*} )$ 
}
is a $\varphi$-fitting pair of integers, then 
$\varphi\sp{\bullet}(\ell\sp{*})= \ell\sp{*}k\sp{*} - \varphi (k\sp{*})$.  
The reverse implication holds for an arbitrary integer-valued function $\varphi$
on ${\bf Z}$.

\begin{proposition} \label{rev-fitting} 
Let $\varphi$ be an arbitrary integer-valued function on ${\bf Z}$, 
and $k\sp{*},\ell\sp{*}$ integers
for which 
$\varphi\sp{\bullet}(\ell\sp{*})= \ell\sp{*}k\sp{*} - \varphi (k\sp{*})$.  
Then the pair 
\REV{
$( k\sp{*},\ell\sp{*} )$ 
}
is $\varphi$-fitting.
\end{proposition}

\Proof 
The definition of $\varphi \sp{\bullet}$ implies that 
\[
\ell\sp{*}k\sp{*} -\varphi (k\sp{*}) 
 = \varphi \sp{\bullet} (\ell\sp{*}) 
 \geq \ell\sp{*}(k\sp{*}+1) -\varphi (k\sp{*}+1),
\] 
from which 
$\varphi (k\sp{*}+1) - \varphi (k\sp{*}) \geq \ell\sp{*}$.  
Analogously, we have
\[
\ell\sp{*}k\sp{*} -\varphi (k\sp{*}) = \varphi \sp{\bullet} (\ell\sp{*}) \geq
\ell\sp{*}(k\sp{*}-1) -\varphi (k\sp{*}-1),
\] 
from which $\ell\sp{*} \geq \varphi (k\sp{*}) - \varphi (k\sp{*}-1)$.  
\finbox \medskip

Proposition~\ref{rev-fitting} results in the following estimation
(that may be viewed as a discrete counterpart of a standard lower
bound in continuous optimization).  

\medskip

\begin{proposition} \label{trivirany} 
Let $R=\{x:  Qx\geq p\}\subseteq {\bf R}\sp{S}$ 
be an integral polyhedron and $\varphi_{s}$ 
an arbitrary integer-valued function on ${\bf Z}$ for each $s\in S$.  
Let $\varphi\sp{\bullet}_{s}$ denote the discrete conjugate of $\varphi_{s}$.  
For any element $z$ of $\odotZ{R}$ and for any integer-valued vector 
$y\geq 0$ (whose components correspond to the rows of $Q$) one has:
\begin{equation} 
 \Phi(z) \geq yp - \Phi\sp{\bullet}(yQ). 
 \label{(trivirany)} 
\end{equation}
\noindent 
If equality holds for $z\sp{*}$ and $y\sp{*}$, then $z\sp{*}$
is a $\Phi$-minimizer of $\odotZ{R}$ and the pair 
\REV{
$( z\sp{*},y\sp{*} )$ 
}
is $\Phi$-compatible.  
\end{proposition}

\Proof 
Let $w:=yQ$.  
By the definition of discrete conjugate, we have
$\varphi \sp{\bullet}_{s}(w(s)) + \varphi_{s}(z(s)) \geq w(s)z(s)$
from which
\begin{equation}
  \Phi(z) = \sum_{s\in S} \varphi_{s}(z(s)) 
  =wz - \big[ \sum_{s\in S} w(s)z(s) - \sum_{s\in S}\varphi_{s}(z(s)) \big] 
   \geq yp - \Phi\sp{\bullet}(yQ).
\label{(estim)} 
\end{equation}

To see the second part, observe that \eqref{(trivirany)} implies that
$\Phi(z) \geq  y\sp{*}p - \Phi\sp{\bullet}(y\sp{*}Q) = \Phi(z\sp{*})$, 
showing that
$z\sp{*}$ is indeed a $\Phi$-minimizer element of $\odotZ{R}$.  
Since we have equality in \eqref{(estim)} for $z\sp{*}$ and $y\sp{*}$, 
it follows for each $s\in S$ that 
$w\sp{*}(s)z\sp{*}(s) - \varphi_{s}(z\sp{*}(s)) 
= \varphi_{s}\sp{\bullet}(w\sp{*}(s))$ where $w\sp{*}:=y\sp{*}Q$.  
By applying
Proposition~\ref{rev-fitting} to $\varphi :=\varphi_{s}$, $\ell\sp{*}:= w\sp{*}(s)$, 
and $k\sp{*}:= z\sp{*}(s)$, we obtain that
\begin{equation}  
\varphi (k\sp{*}) - \varphi (k\sp{*}-1) 
\leq \ell\sp{*} \leq \varphi (k\sp{*}+1) - \varphi (k\sp{*}), 
\end{equation}
that is,
\begin{equation}
 \varphi_{s}(z\sp{*}(s)) - \varphi_{s}(z\sp{*}(s)-1)
 \leq w\sp{*}(s)
 \leq \varphi_{s}(z\sp{*}(s)+1) - \varphi_{s}(z\sp{*}(s)), 
\end{equation}
and hence the pair 
\REV{
$( z\sp{*},w\sp{*} )$ 
}
is $\Phi$-fitting,
showing that 
\REV{
$( z\sp{*},y\sp{*} )$ 
}
is $\Phi$-compatible.  
\finbox

\begin{theorem} \label{conjmm} 
Under the same assumptions as in Theorem~{\rm \ref{min-max-uj}}, 
one has the following min-max formula:
\begin{equation}  
\min \{ \Phi (z) :  z\in \odotZ{R} \} 
= \max \{ yp - \Phi\sp{\bullet}(yQ):  y\geq 0 \ \ \hbox{\rm integer-valued} \}.  
\label{(dualmaxD)} 
\end{equation}
\noindent 
The optimal dual vector $y$ can be chosen so as to have at most $2\vert S\vert$
positive components.  
\end{theorem}

\Proof 
\REV{
Let $z\sp{*}\in \odotZ{R}$ be a minimizer element in \eqref{(dualmaxuj)}.
Let $y\sp{*}$ be a non-negative integer-valued vector, guaranteed in Corollary~\ref{Optcrit}, 
that meets the optimality criteria in \eqref{(optcrit1bD)} and \eqref{(optcrit2bD)}.  
By \eqref{(optcrit2bD)}, $( z\sp{*},w\sp{*} )$ with $w\sp{*}:=y\sp{*}Q$
is a $\Phi$-fitting pair.
By \eqref{(optcrit1bD)} we have
}
\begin{equation} 
 \Phi(z\sp{*})= y\sp{*}p  - [ (y\sp{*}Q)z\sp{*} - \Phi(z\sp{*})].
\label{(egyen)} 
\end{equation}

For $s\in S$, consider $\varphi_{s}$ and its discrete conjugate
$\varphi \sp{\bullet}_{s}$.  
For $k\sp{*}:=z\sp{*}(s)$ and $\ell \sp{*}:=w\sp{*}(s)$,
\REV{
$( k\sp{*},\ell\sp{*} )$ 
}
is a $\varphi_{s}$-fitting pair.  
Proposition~\ref{fitting}, when applied to $\varphi_{s}$ in place of $\varphi$,
\REV{
shows that
}
\begin{equation} 
 \ell \sp{*} k\sp{*} - \varphi_{s}(k\sp{*}) 
  = \max \{ \ell \sp{*}k - \varphi_{s}(k) :  k\in {\bf Z}\} 
  = \varphi \sp{\bullet}_{s}(\ell\sp{*}), 
\end{equation}
from which 
$(y\sp{*}Q) z\sp{*} - \Phi (z\sp{*}) =
\REV{
\Phi\sp{\bullet} (w\sp{*}) 
}
= \Phi\sp{\bullet} (y\sp{*}Q)$ 
follows.  
\REV{
By substituting this into \eqref{(egyen)}
we obtain \eqref{(dualmaxD)}.
}
\finbox

\begin{corollary} \label{ketfele} 
Let $z\sp{*}$ be 
\RED{a}
$\Phi$-minimizer element of $\odotZ{R}$.  
If 
\REV{
$(y\sp{*}, z\sp{*})$ is an optimal solution to \eqref{(dualmaxuj)},
}
then $y\sp{*}$ is an optimal solution to \eqref{(dualmaxD)}.  
If $y\sp{*}$ is an optimal solution to \eqref{(dualmaxD)}, 
then the pair 
\REV{
$( z\sp{*},y\sp{*} )$ 
}
is $\Phi$-compatible 
and 
$(y\sp{*}, z\sp{*})$ 
 is an optimal solution to \eqref{(dualmaxuj)}.  
\end{corollary}

\Proof 
The first part is an immediate consequence of the proof of Theorem~\ref{conjmm}.  
The second part follows from Theorem~\ref{conjmm} and the second half of Proposition~\ref{trivirany}.
\finbox

\begin{remark} \rm \label{Philinear.3} 
In the special case when $\Phi$ is linear and defined by $\Phi(w)=cw$, 
one can easily observe that
$\Phi\sp{\bullet}(w)=0$ when $w=c$ and $\Phi\sp{\bullet}(w)$ 
has a $+\infty $ summand when $w\not =c$.  
Therefore $\Phi\sp{\bullet}(yQ)$ in \eqref{(dualmaxD)} is
finite only if $yQ=c$ and in this case $\Phi\sp{\bullet}(yQ)=0$.  
This means that the maximum in \eqref{(dualmaxD)} is equal to 
$\max \{ yp :  yQ=c, y\geq 0\}$, 
showing that Theorem~\ref{conjmm} also specializes to the
integral version of the linear programming duality theorem formulated
for box-TDI polyhedra.  
\hfill $\bullet$ 
\end{remark}

The results above can be extended to the case when $R$ is defined 
by a box-TDI system 
$[ Q'x\geq p', \  Q\sp{=}x=p\sp{=} ]$
 because
\REV{
this means, by definition, 
that the system 
$[ Q'x\geq p', \ Q\sp{=}x\geq p\sp{=},  \  -Q\sp{=}x\geq -p\sp{=} ]$ 
is also box-TDI and defines the same polyhedron $R$.  
}%
We call a dual vector $y=(y',y\sp{=})$ {\bf sign-feasible} 
if $y'\geq 0$. 
That is, we require non-negativity of
those components that correspond to the rows of $Q'$.

\begin{theorem} \label{min-max-gen} 
Suppose that in Theorem {\rm \ref{min-max-uj}}
 the box-TDI polyhedron is given in form 
$R=\{x:  Q'x \geq p', Q\sp{=}x =p\sp{=}\}$, 
where each of  $Q'$, $Q\sp{=}$, $p'$, $p\sp{=}$ is
integer-valued.  Then 
\begin{align}
& \min \{ \Phi (z) :  z\in \odotZ{R} \}
\nonumber 
\\ &
=   \max \{ \Phi(z) - y(Qz-p) :  \ z\in \odotZ{R} , \ y \ \ 
\hbox{\rm sign-feasible and integer-valued, 
\REV{
$( z,y )$ 
}
$\Phi$-compatible} \} 
\nonumber 
\\ &
=  \max \{ yp - \Phi\sp{\bullet}(yQ):  y \ \ \hbox{\rm sign-feasible and integer-valued} \}, 
\label{(dualmaxD-gen)} 
\end{align}
where 
$Q =\begin{pmatrix}Q'\\Q\sp{=}\end{pmatrix}$ 
and $p =\begin{pmatrix}p'\\p\sp{=}\end{pmatrix}$.  
An element 
$z\sp{*}\in \odotZ{R}$ is a $\Phi$-minimizer if and only if there exists a
sign-feasible integer-valued vector $y\sp{*}$ meeting the optimality criteria:
\begin{align}  
& y\sp{*}(Qz\sp{*}-p) = 0, 
\label{(optcrit1bD2)} 
\\ &
 \Phi' (z\sp{*}- {\bf \underline{1}} ) \leq y\sp{*}Q \leq \Phi' (z\sp{*}).
\label{(optcrit2bD2)} 
\end{align}
\noindent 
Moreover, $y\sp{*}$ can be chosen in such a way that the number 
of its non-zero components is at most 
$2\vert S\vert$.
\finbox 
\end{theorem}

We formulate yet another variant for the maximum in the min-max theorem.  
This version is useful in cases when there is a simple formula for 
\REV{
$\mu_{R}(w):= \min \{wx:  x\in R\}$, 
}
see the next section on special box-TDI polyhedra.

\begin{theorem} \label{muR} 
Let $R=\{x:  Q' x \geq p', Q\sp{=} x = p\sp{=} \}$ 
be a box-TDI polyhedron, where each of $Q'$, $Q\sp{=}$, $p'$, $p\sp{=}$ is
integer-valued.  Let $\Phi$ be an integer-valued separable discrete
convex function which is bounded from below on $R$.  Then
\begin{equation} 
 \min \{ \Phi (z) :  z\in \odotZ{R} \} 
= \max \{\mu_{R}(w) - \Phi\sp{\bullet}(w):  w\in {\bf Z} \sp{S}\}.  
\label{(minmax-muR)} 
\end{equation} 
\end{theorem}

\Proof 
For $z\in \odotZ{R}$ and $w\in {\bf Z} \sp{S}$, we have
\begin{equation}  
\Phi(z) = wz - [wz- \Phi(z)] \geq \mu_{R}(w) - \Phi\sp{\bullet}(w),
\label{(muR-estim)} 
\end{equation} 
from which $\min \geq \max $ follows.

To see the reverse direction, we show that 
there is an element $z\sp{*}$ of $\odotZ{R}$ and an integral vector $w\sp{*}$ meeting 
\eqref{(muR-estim)} with equality.  
Let $z\sp{*}$ be a $\Phi$-minimizer of $\odotZ{R}$, 
$y\sp{*}$ a maximizer in \eqref{(dualmaxD-gen)}, and let
$w\sp{*}:=y\sp{*}Q$.  
Then $\Phi(z\sp{*}) = py\sp{*} - \Phi\sp{\bullet}(y\sp{*}Q)$
holds by Theorem~\ref{min-max-gen}, and a straightforward estimation
(the weak duality theorem of linear programming) shows that 
$\mu_{R}(w\sp{*}) \geq y\sp{*}p$.  
This and \eqref{(muR-estim)} 
(when applied to $z\sp{*}$ and $w\sp{*}$) imply
\begin{equation}  
 \Phi(z\sp{*}) 
 \geq \mu_{R}(w\sp{*}) - \Phi\sp{\bullet}(w\sp{*}) 
 \geq y\sp{*}p  - \Phi\sp{\bullet}(y\sp{*}Q) 
 = \Phi(z\sp{*}), 
\end{equation}
from which equality follows throughout, and hence \eqref{(minmax-muR)} holds indeed.  
\finbox

\begin{remark} \rm \label{most-general} 
Theorem~\ref{min-max-gen} can be further extended 
to the formally more general framework 
where primal non-negativity constraints are written separately.  
In this case, the primal polyhedron $R$ is defined by a box-TDI system as follows:
\[
 R:=\{(x_{1},x_{2}):  
 Q'x_{1}+ A'x_{2} \geq p', \ 
 Q\sp{=} x_{1} + A\sp{=}x_{2} = p\sp{=}, \ 
 x_{2}\geq 0\}.
\]  
The min-max theorem for the minimum of $\Phi$ and the
optimality criteria, though technically more complex, can also be
described by applying Theorem~\ref{min-max-gen}.  
\hfill $\bullet$ 
\end{remark}


\section{Special box-TDI polyhedra} 
\label{spec.poly}

In this section, we consider special box-TDI polyhedra.

\subsection{Polyhedra defined by TU-matrices}

It is known that if 
$Q =\begin{pmatrix}Q'\\Q\sp{=}\end{pmatrix}$ 
is a totally unimodular (TU) matrix and 
$p =\begin{pmatrix} p' \\ p\sp{=} \end{pmatrix}$ 
is an integral vector, then the linear system
$[ Q'x\geq p', \ Q\sp{=}x=p\sp{=} ]$ 
(and the polyhedron 
$\{x:  Q'x\geq p', \  Q\sp{=}x=p\sp{=}\}$) is box-TDI.  
But a TU-matrix $Q$ may define box-TDI polyhedra in other ways, as well.

\begin{proposition} \label{boximage} 
Let $Q$ be a totally unimodular matrix,
and let $f\leq g$ be integer-valued bounding vectors (of appropriate dimension) 
where $f$ may have $-\infty $ while $g$ may have $+\infty $ components.  
Then the polyhedron 
\[
R':=\{z:  z = Qx \ \ \mbox{\rm for some $x$ meeting} \  \  f \leq x \leq g\}
\]
 is box-TDI.  
Analogously, if $\ell\leq u$
are integer-valued bounding vectors (of appropriate dimension), 
then the polyhedron
\[
R'':=\{w:  w = yQ \ \ \mbox{\rm for some $y$ meeting} \ \ \ell \leq y \leq u\}
\]
is box-TDI.  
\end{proposition}

\Proof 
Since the operation of adding a unit vector $(1,0,\dots ,0)$ to $Q$ 
as a new row or a new column preserves total unimodularity, 
the system  $[  Qx -z =0, \ f\leq x\leq g ]$ is box-TDI.  
But then $R'$ is the projection of the polyhedron 
$\{(x,z):  Qx -z =0, \ f \leq x\leq g\}$ 
along the coordinate axes of $x$ (to the components of $z$), 
and since projection, by Proposition~\ref{projection}, 
preserves box-TDI-ness, $R'$ is indeed box-TDI.  The second part follows from
the first one since the transpose of a 
\REV{
TU-matrix
}
 is also totally unimodular.  
\finbox

\medskip

It follows that Theorem~\ref{muR} (for example) can be applied 
to the box-TDI polyhedra occurring in Proposition~\ref{boximage}.  A special
case is the polyhedron of feasible flows defined on the edge-set of a
digraph $D=(V,A)$ by 
$\{x\in {\bf R} \sp A:  \varrho _x(v)-\delta_x(v) = m(v)$
for each $v\in V, \ f\leq x\leq g\}$,
where
$m:V\rightarrow {\bf Z}$ is a function on $V$ with $\widetilde m(V)=0$ 
while $f:A\rightarrow {\bf Z} \cup \{-\infty \}$ and
$g:A\rightarrow {\bf Z} \cup \{+\infty \}$ 
are bounding functions on $A$ with $f\leq g$.  
Here $\varrho _x(v):=\sum [x(uv):  uv\in A]$ and
$\delta _x(v):=\sum [x(vu):  vu\in A]$.  
The classic notions of feasible $st$-flows with given flow-amount 
as well as feasible circulations fit into this framework.  
Another special case of TU-polyhedra is the one of feasible potentials.

A network matrix $Q$ is a more general TU-matrix which is defined by a
digraph whose underlying undirected graph is connected.  
For a spanning tree 
$T$ of $D$, the rows of $Q$ correspond to the elements
of $T$, the columns correspond to the edges in $A-T$,
and the column corresponding to $e$ is the signed characteristic vector of the
fundamental circuit belonging to $e$.

\subsection{M-convex and M$_{2}$-convex sets}
\label{SCmsetm2set}

Let $B:=B'(p)=\{\widetilde x(Z)\geq p(Z)$ for $Z\subset S$ and
$\widetilde x(S)=p(S)\}$ be 
\RED{
the
}
base-polyhedron defined by an integral
supermodular function $p$ for which $p(S)$ is finite.  Recall that the
set $\odotZ{B}$ of integral elements of $B$ 
is called an M-convex set.
A basic property of base-polyhedra is that they are box-TDI.

Recall that the linear extension (Lov\'asz extension) $\hat p$ of $p$ is defined by
\begin{equation}  \label{lovextdef}
\hat p(w) := p(S_n)w(s_n) + \sum_{j=1}\sp {n-1} p(S_{j})[w(s_{j}) - w(s_{j+1})],
\end{equation}
where $n=\vert S\vert $, the elements of $S$ are indexed in
such a way that 
$w(s_{1})\geq \cdots \geq w(s_n)$, and 
$S_{j}=\{s_{1},\dots ,s_{j}\}$ for $j=1,\dots ,n$.  
(Here $p(S_{j})[w(s_{j}) - w(s_{j+1})]$ is defined $0$ when 
$w(s_{j}) - w(s_{j+1})=0$ even if $p(S_{j})$ is not finite.)

Recall the definitions of $z$-tight sets and strict $w$-top sets.  
For
\REV{
a
}
supermodular function $p$, a theorem of Edmonds \cite{Edmonds70} is as follows.

\begin{claim} 
\label{greedy} 
Let $p:2\sp{S}\rightarrow {\bf Z} \cup \{-\infty \}$ 
be a supermodular function for which $p(S)$ is finite.  
For an integral cost-function $w$ on $S$, one has 
\[ 
\hat p(w) = \mu _B(w)= \mu_{\odotZ{B}}(w),
\] 
where $\mu _B(w):=\min \{wx:  x\in B\}$ and 
$\mu_{\odotZ{B}}(w):=\min \{wx:  x\in \odotZ{B}\}$.  
In particular, 
$\hat p(w)= -\infty $ if and only if $wz$ is unbounded from below over $\odotZ{B}$.  
When $\hat p(w)>-\infty $, an element $z\in \odotZ{B}$ is
a $w$-minimizer if and only if each strict $w$-top set is $z$-tight.
\finbox 
\end{claim}

By combining Claim \ref{greedy} with Theorem~\ref{muR}, we arrive at
the starting min-max formula described in Theorem~\ref{min-max.orig}.

Theorem~\ref{min-max-M2} can also be derived in an analogous way.  
Let $B_{1}:=B'(p_{1})$ and $B_{2}:=B'(p_{2})$ 
be base-polyhedra defined by integer-valued supermodular functions 
$p_{1}$ and $p_{2}$ for which $B:=B_{1}\cap B_{2}$ is non-empty.  
A fundamental theorem of Edmonds states that $B$ is box-TDI.  
A version of the well-known weight-splitting theorem 
(\cite{Frank-book}, Theorem 16.1.8) states
for an integral vector $w$ that
\[
 \mu _B(w) = \mu_{\odotZ{B}}(w) 
 = \max \{\hat p_{1}(w_{1}) + \hat p_{2}(w_{2}) :  
  w_{1}+w_{2}=w, \ w_{1}, w_{2} \hbox{ integral} \}.  
\]
\noindent 
Combining this formula with Theorem~\ref{muR}, we arrive at
Theorem~\ref{min-max-M2}.

It should be noted that the intersection of two integral g-polymatroids 
is also box-TDI and so is a submodular flow polyhedron
(by a theorem of Edmonds and Giles \cite{Edmonds-Giles}).  
Therefore the general min-max formulas described in Theorem~\ref{conjmm} 
can be specialized to these cases as well.

\subsection{Direct proof for M-convex sets} \label{min-max.proof}

The goal of this section is to provide a direct proof of the
non-trivial part of Theorem~\ref{min-max.orig}.  
The proof is independent of the results in Section \ref{main-min-max} 
and gives rise to a strongly polynomial algorithm to compute the optimal dual,
provided that an optimal solution to the primal problem is available.
Namely, we prove the following.

\begin{theorem} \label{nontrivM} 
Let $\Phi$ be an integer-valued separable discrete convex function.  
Let $z\sp{*}$ be a $\Phi$-minimizer element of an M-convex set $\odotZ{B}$ 
defined by a finite-valued supermodular function $p$.  
There exists an integer-valued vector $w\sp{*}\in {\bf Z}\sp{S}$ 
for which $z\sp{*}$ and $w\sp{*}$ meet the optimality criteria
\eqref{(optcrit1b)} and \eqref{(optcrit2d)} 
(or \eqref{(optcrit2.conc}).  
\end{theorem}

\Proof 
\REV{
First recall that the $z\sp{*}$-tight sets form a ring-family (lattice)
closed under intersection and union.
}
For $s\in S$, let $T(s)$ denote the unique smallest $z\sp{*}$-tight set containing $s$.

\begin{claim} \label{Goenb} 
For a $\Phi$-minimizer element $z\sp{*}$, 
\begin{equation}  
 \varphi_{s}'(z\sp{*}(s)-1) \leq \varphi_{t}'(z\sp{*}(t)) 
\label{(Groen1b)}
\end{equation}
holds whenever $t\in T(s)$.  
\end{claim}

\Proof 
\RED{
First we show that $z':=z\sp{*} - \chi_s + \chi_t$ belongs to $\odotZ{B}$, 
which is equivalent to requiring that 
$\widetilde z'(Z)\geq p(Z)$ 
for every subset $Z\subseteq S$.  
Indeed, if $Z$ is not $z\sp{*}$-tight, then 
$\widetilde z'(Z)\geq \widetilde z\sp{*}(Z)-1\geq p(Z)$; 
if $Z$ is $z\sp{*}$-tight and $s\not \in Z$, then
$\widetilde z'(Z)\geq \widetilde z\sp{*}(Z)=p(Z)$.  
Finally, if $Z$ is $z\sp{*}$-tight and $s\in Z$, 
then $T(s)\subseteq Z$, and hence $s,t\in Z$, from which 
$\widetilde z'(Z)=\widetilde z\sp{*}(Z)=p(Z)$.
}

\RED{
As $z'$ belongs to $\odotZ{B}$, we have
$\Phi(z\sp{*})\leq \Phi(z')$, from which
}
\[ 
\varphi_{s}(z\sp{*}(s)) + \varphi_{t}(z\sp{*}(t)) 
\leq \varphi_{s}(z\sp{*}(s)-1) + \varphi_{t}(z\sp{*}(t)+1), 
\]
that is,
\begin{equation}  
\varphi_{s}(z\sp{*}(s)) - \varphi_{s}(z\sp{*}(s)-1) 
\leq \varphi_{t}(z\sp{*}(t)+1) - \varphi_{t}(z\sp{*}(t)), 
\label{(Groen1)} 
\end{equation}
which is exactly \eqref{(Groen1b)}.

By the discrete convexity of $\varphi_{s}$, we have 
$\varphi_{s}(z\sp{*}(s)) - \varphi_{s}(z\sp{*}(s)-1) 
  \leq \varphi_{s}(z\sp{*}(s)+1) - \varphi_{s}(z\sp{*}(s))$, 
implying \eqref{(Groen1)} (and hence \eqref{(Groen1b)})
for the case $s=t$.  
\finbox

\medskip

Our goal is to find an integer-valued $w\sp{*}$ meeting the optimality
criteria in the theorem.  
Define $w\sp{*}$ as follows:
\begin{equation}  
 w\sp{*}(s):= \min \{\varphi_{t}'(z\sp{*}(t)):  t \in T(s)\}.
\label{(pidef)} 
\end{equation}

\begin{claim}  
 \label{optcrit2}$w\sp{*}$ and $z\sp{*}$ meet 
\REV{
the optimality
}
criterion \eqref{(optcrit2d)}.  
\end{claim}

\Proof 
The definition of $w\sp{*}(s)$ in \eqref{(pidef)} and $s\in T(s)$ imply that 
$w\sp{*}(s)= \min \{\varphi_{t}'(z\sp{*}(t)):  t \in T(s)\} 
\leq \varphi_{s}'(z\sp{*}(s))$, 
from which $w\sp{*}(s) \leq \varphi_{s}'(z\sp{*}(s))$ 
follows.  Furthermore, \eqref{(pidef)} and \eqref{(Groen1b)} imply that 
$w\sp{*}(s)= \min \{\varphi_{t}'(z\sp{*}(t)): t\in T(s)\} \geq \varphi_{s}'(z\sp{*}(s)-1)$.  
Hence \eqref{(optcrit2d)} holds.  
\finbox

\begin{claim}  
\label{optcrit1}$w\sp{*}$ and $z\sp{*}$ meet 
\REV{
the optimality
}
 criterion \eqref{(optcrit1b)}.  
\end{claim}

\Proof 
Let $\beta_{1}>\beta_{1}>\cdots >\beta _\ell$ denote 
the distinct values of the components of $w\sp{*}$, 
and let $C_{i}:= \{s:  w\sp{*}(s)\geq \beta_{i}\}$ for 
\REV{
$i=1,\dots ,\ell$.  
}
Let $S'_{1}:=C_{1}$ and
$S'_{i}:=C_{i}-C_{i-1}$ for 
$i=2,\dots ,\ell$.  
Then ($\emptyset \not =$)
\ $C_{1}\subset C_{2}\subset \cdots \subset C_\ell$ ($=S$) 
is a chain whose members are the strict $w\sp{*}$-top sets, 
while $\{S'_{1},\dots ,S'_\ell\}$ 
is a partition of $S$ for which $w\sp{*}(s)=\beta_{i}$ holds for every $s\in S'_{i}$.

For every $t\in T(s)$, we have $T(t)\subseteq T(s)$ and hence 
$w\sp{*}(t) \geq w\sp{*}(s)$, 
implying that $T(s)\subseteq C_{i}$ whenever $s\in S'_{i}$.  
Therefore 
\REV{
$C_{i}=\bigcup_{s \in S} \{T(s):  w\sp{*}(s)\geq \beta_{i}\}$ 
}
and hence each $C_{i}$ is $z\sp{*}$-tight, 
showing that 
\REV{
the optimality
}
criterion \eqref{(optcrit1b)} holds.  
\finbox

\medskip

As $z\sp{*}$ and $w\sp{*}$ meet the optimality criteria, the proof of
the theorem is complete.  
\BB

\medskip 

In order to compute $w\sp{*}$, we have to be able to determine
the unique smallest $z\sp{*}$-tight set $T(s)$ containing an element $s\in S$.  
This is easy once we are able to decide for a given pair
$\{s,t\}$ of elements of $S$ whether there is a $z\sp{*}$-tight  $s\overline{t}$-set.  
But this can be done by minimizing the submodular function
$\widetilde z\sp{*}-p$ over the $s\overline{t}$-sets, 
which is doable in strongly polynomial time with the help 
of a general subroutine to minimize a submodular function
\REV{
\cite{McC06handbk}.
}

\begin{remark} \rm 
We formulated and proved Theorem~\ref{nontrivM} 
for the special case when the defining supermodular function $p$ is finite-valued.  
But the arguments above can easily be extended to the general case 
when $p$ may have $-\infty $ values (but preserving the finiteness of $p(S)$), 
that is, $B'(p)$ may be unbounded.  
\hfill $\bullet$
\end{remark}


\section{Special discrete convex functions}
\label{SCspedcfn}

\subsection{Minimizing the square-sum}
\label{SCsqsum}

Consider the special case when 
$\varphi_{s}(k) := \varphi (k):= k\sp{2}$
for each $s\in S \ (= \{1,2,\dots ,n\})$ 
and hence the separable discrete convex function $\Phi$ to be minimized 
is given by
$\Phi(z):=z\sp{2}$, where $z\sp{2}=\sum [z(i)\sp{2}:  i=1,2,\dots ,n]$.
That is, we want to minimize the square-sum of the components of $z$.
For this problem Theorem~\ref{min-max-gen} 
\REV{
is specialized
}
as follows.

The discrete conjugate of $\varphi (k)= k\sp{2}$
is explicitly available, namely,
for integer $\ell$:
\[
  \varphi \sp{\bullet}(\ell) = 
  \llfloor  \frac{\ell}{2}\rrfloor  
  \  \llceil \frac{\ell}{2}\rrceil 
\]
(the proof is immediate and is given explicitly in a more general case in Proposition \ref{PRakk.conj} in Section~\ref{SCwtdsqsum}).
Then, for an integral vector $w\in {\bf Z}\sp{S}$, we obtain
$\Phi\sp{\bullet}(w) = \lfloor  w / 2 \rfloor  \  \lceil  w / 2  \rceil$,
where, for a vector  $x=(x(1),\dots ,x(n))$,
we use notations  
$\lfloor  x \rfloor := (\lfloor  x(1) \rfloor  , \dots , \lfloor  x(n) \rfloor  )$ 
and
$\  \lceil x \rceil := (\lceil x(1) \rceil , \dots ,  \lceil x(n) \rceil)$.  
In addition, 
we observe for the condition \eqref{(optcrit2bD2)} that 
$\varphi '(k)= (k+1)\sp{2} - k\sp{2}= 2k+1$ and
$\varphi '(k-1)= 2k-1$.  
Hence, for an integral vector $z\in {\bf Z}\sp{S}$, 
we have $\Phi'(z) = 2z+{\bf \underline{1}}$ and
$\Phi'(z- {\bf \underline{1}} ) = 2z -{\bf \underline{1}}$,
where ${\bf \underline{1}} = \chi_{S}$.

\begin{theorem} \label{min-maxW} 
Let $Q =\begin{pmatrix}Q'\\Q\sp
=\end{pmatrix}$ be an integral matrix and 
$p = \begin{pmatrix}p' \\ p\sp{=} \end{pmatrix}$ 
an integral vector, and suppose that the linear system
$[ Q'x\geq p', \ Q\sp{=}x = p\sp{=} ]$ is box-TDI.  
Let $R:=\{x:  Q'x\geq p', Q\sp{=}x = p\sp{=} \} \subseteq {\bf R}\sp{S}$ 
be the (box-TDI) polyhedron defined by this system.  
Then 
\begin{align} 
 & \min \{ z\sp{2} :  z\in \odotZ{R} \} 
\label{(Q1Q2mini)} 
\\ &
= \max \{ yp - 
   \llfloor  \frac{yQ}{2} \rrfloor  \  \llceil \frac{yQ}{2} \rrceil
   :  \ y=(y',y\sp{=}) \ \  \hbox{\rm sign-feasible and integer-valued}\},
\label{(dualmaxW)} 
\end{align} 
where the sign-feasibility of $y$ means that $y' \geq 0$.  
Moreover, an integral element $z\sp{*}\in \odotZ{R}$ 
is a square-sum minimizer 
if and only if 
there exists a sign-feasible integral vector $y\sp{*}$ 
for which the following optimality criteria hold:
\begin{align} 
& y\sp{*} (Qz\sp{*}-p) =0, 
\label{(optcrit1bW} 
\\ &
  2z\sp{*} - {\bf \underline{1}} \ \leq \ y\sp{*}Q \leq \ 2z\sp{*} + {\bf \underline{1}}.
\label{(optcrit2dW)} 
\end{align} 
\noindent 
The optimal (integral) dual solution $y\sp{*}$ can be chosen
in such a way that the number of its non-zero components is at most 
$2\vert S\vert$.
\finbox 
\end{theorem}

\begin{remark} \rm \label{sqsumsubgrcond} 
It is worth noting that \eqref{(optcrit2dW)} is equivalent to
\begin{equation} 
 \llfloor  \frac{y\sp{*}Q}{2}\rrfloor     
 \ \leq \ z\sp{*} \ \leq \ 
 \  \llceil \frac{y\sp{*}Q}{2}\rrceil . 
\label{(ekvopt)} 
\end{equation}
\hfill $\bullet$ 
\end{remark}

\begin{remark} \rm \label{trivest} 
For a simple understanding, it is worth providing a direct proof 
of the trivial inequality $\min\geq \max$
that relies neither on $\Phi$-compatibility nor on conjugacy.  
For real vectors $w$ and $z$ in ${\bf R}\sp{n}$, one has the obvious
estimation $z(w-z) \leq (w/2)(w/2)$.  
For integral vectors $w$ and $z$, 
the stronger inequality 
$z(w-z) \leq \lfloor  w / 2 \rfloor  \,  \lceil w / 2 \rceil$ 
holds, with equality precisely if 
$\lfloor w/2 \rfloor  \leq z \leq   \lceil w/ 2 \rceil$, 
that is,
$z(s)\in \{ \lfloor  w(s) / 2  \rfloor  ,  \lceil w(s) / 2 \rceil \}$ 
for each $s\in S$.  
This implies for any $z\in \odotZ{R}$
and for any integral vector $y=(y',y\sp{=})$ with $y'\geq 0$ that
\begin{equation}  
z\sp{2} \ = \ (yQ)z - ((yQ)z-z\sp{2}) 
   \ = \ y(Qz) - (yQ-z)z 
   \ \geq \  yp - \llfloor  \frac{yQ}{2}\rrfloor  \  \llceil \frac{yQ}{2}\rrceil ,
\label{(Restimb)} 
\end{equation}
from which $\min\geq \max$ follows.  Moreover,
equality holds in \eqref{(Restimb)} for $z\sp{*}$ and $y\sp{*}$ in place
of $z$ and $y$ precisely if the optimality criteria \eqref{(optcrit1bW}
and \eqref{(ekvopt)} hold.  
\hfill $\bullet$ 
\end{remark}

\begin{remark} \rm \label{affine} 
Minimizing the square-sum over an affine subspace $R=\{x:  Qx=p\}$ 
is a standard problem of linear algebra.
If $R$ is an integral box-TDI affine subspace 
(that is, $Q'$ is empty in Theorem~\ref{min-maxW}), 
then we get the following min-max formula:
\begin{equation}  
\min \{ z\sp{2} :  z\in \odotZ{R} \} 
 = \max \{ yp - 
 \llfloor  \frac{yQ}{2}\rrfloor  \  \llceil \frac{yQ}{2} \rrceil : 
 \ y \ \hbox{\rm integer-valued} \}.  
\label{(affin.min.max)} 
\end{equation}
\noindent 
For the case when $Q$ is totally unimodular, McCormick et al.~\cite{MPSV20} 
described a polynomial algorithm for computing the minimum.  
\hfill $\bullet$ 
\end{remark}

\begin{remark} \rm \label{Sprime} 
Theorem~\ref{min-maxW} can easily be extended 
to the slightly more general case when the goal is 
to minimize the sum of squares over a given subset $S'$ of coordinates.
In this case \eqref{(Q1Q2mini)} turns to
\[
 \min \{ \sum_{s\in S'} z(s)\sp{2}:  z\in \odotZ{R} \}. 
\]
\noindent 
Let $Q''$ denote a matrix consisting of the columns of
$Q=\begin{pmatrix} Q' \\ Q\sp{=} \end{pmatrix}$ 
corresponding to the elements of $S-S'$.  
Then \eqref{(dualmaxW)} transforms to the following:
\begin{equation} 
 \max \{ yp - 
  \llfloor  \frac{yQ}{2}\rrfloor  \  \llceil \frac{yQ}{2}\rrceil
:  \ y \ \ \hbox{\rm sign-feasible and integer-valued}, \  yQ''=0\} .
\end{equation}
\hfill $\bullet$ 
\end{remark}

\subsection{Flows and circulations}
\label{SCflowcirc}

In this section, we specialize Theorem~\ref{min-maxW} to network flows.  
Let $D=(V,A)$ be a digraph and let $m$ be an integral function on $V$ 
for which $\widetilde m(V)=0$.  
A function $x$ on $A$ is called an {\bf $m$-flow} if 
\begin{equation}  
\varrho _x(v) - \delta _x(v) = m(v) \ \ \hbox{\rm for every}\ \ v\in V. 
\label{(mflow)} 
\end{equation}
\noindent 
Note that this is equivalent to $Q_Dx=m$
where $Q_D$ denotes the signed incidence matrix of $D$.  
The columns of $Q_D$ correspond to the edges of $D$ 
while the rows correspond to the nodes. 
\RED{
An entry of $Q_D$, corresponding to edge $a$ and node $v$,
}
is $+1$ or $-1$ according as $a$ enters or leaves $v$, and $0$ otherwise.  
By the assumption $\widetilde m(V)=0$, \eqref{(mflow)} is equivalent to
\begin{equation}
  \varrho _x(v) - \delta _x(v) \geq m(v) \ \ \hbox{\rm for every}\ \ v\in V, 
\label{(mflow2)} 
\end{equation}
or concisely $Q_Dx\geq m$.

By Hoffman's circulation theorem, there is a non-negative integral $m$-flow 
if and only if 
$\widetilde m(X)\geq 0$ holds for every subset $X\subseteq V$ for which $\delta _D(X)=0$.  
We assume that there is a non-negative integral $m$-flow $z$ 
and we want to characterize those minimizing the square-sum 
$z\sp{2}=\sum [z(a)\sp{2}:a\in A]$.  
We are going to specialize Theorem~\ref{min-maxW}.  In this case, $y$ is a
$(\vert V\vert +\vert A\vert)$-dimensional vector but in order to have
a better fit to the standard notation in network flow theory, we
replace $y$ by a vector $(\pi ,h)$ where $\pi $ (a \lq potential\rq)
is defined on $V$ while $h$ is defined on $A$.

\begin{theorem} \label{minsq.flow} 
The minimum square-sum of a non-negative integral $m$-flow is equal to
\begin{align} 
&  \max \{m\pi - \llfloor  \frac{\max(\Delta_{\pi} ,0)}{2} \rrfloor 
\  \llceil \frac{\max(\Delta_{\pi} ,0)}{2}\rrceil :  \pi :  V\rightarrow {\bf Z}_{+} \}, 
\label{(mflowdualx)} 
\end{align}
where $\Delta_{\pi} $ denotes the 
\REV{
tension (= potential-difference)
}
defined by $\pi $, that is, $\Delta_{\pi} (uv)= \pi (v)-\pi (u)$ for
every edge $uv\in A$, or concisely, $\Delta_{\pi} =\pi Q_D$.  
The minimum square-sum of an integral $m$-flow is equal to
\begin{equation}  
 \max \{m\pi - \llfloor  \frac{\Delta_{\pi}}{2} \rrfloor  
      \  \llceil \frac{\Delta_{\pi}}{2}\rrceil :  \pi :  V\rightarrow {\bf Z}_{+} \}.
\label{(mflowdual.affine)} 
\end{equation}
\end{theorem}

\Proof 
Apply Theorem~\ref{min-maxW} to the special case 
when the system is $Q'x\geq p'$ (and $Q\sp{=}$ is empty) where
$Q'=\begin{pmatrix} Q_D \\ I \end{pmatrix}$ 
and $p'$ is defined by
$p'(v):=m(v)$ for $v\in V$ and $p'(a):=0$ when $a\in A$.  
(Here $I$ denotes the $\vert A\vert $ by $\vert A\vert $ unit-matrix).  
The optimal dual vector $y=y'$ in Theorem~\ref{min-maxW} can be written 
in the form $y=(\pi ,h)$, where $\pi $ corresponds to the sub-vector of $y$ 
whose components are assigned to the rows of $Q_D$ 
(that is, to the nodes of $D$) while the components of $h$ are assigned to the rows
of $I$ (that is, to the edges of $D$).
Then the expression \eqref{(dualmaxW)} in Theorem~\ref{min-maxW} 
takes the following form
\begin{equation}
 \max \{m\pi -
   \llfloor  \frac{\Delta_{\pi} + h}{2} \rrfloor  
 \ \llceil \frac{\Delta_{\pi} + h}{2}\rrceil : 
   \pi :  V\rightarrow {\bf Z}_{+}, \ h: A\rightarrow {\bf Z}_{+} \} .
\label{(mflowdual)} 
\end{equation}
To see that this is equal to \eqref{(mflowdualx)}, 
it suffices to observe that in an optimal solution $(\pi,h)$ to \eqref{(mflowdual)}, 
if $\Delta_{\pi} (a)$ is negative for an edge $a$ of $D$, 
then $h(a)$ may be chosen to be $\vert \Delta_{\pi} (a)\vert $, 
while if $\Delta_{\pi} (a)$ is non-negative, then $h(a)$
may be chosen to be zero, and hence 
$\Delta_{\pi} + h = \max\{ \Delta_{\pi} ,0\}$.
The expression \eqref{(mflowdual.affine)} follows analogously from
\eqref{(affin.min.max)} in Remark \ref{affine}.  
\finbox

\begin{remark} \rm
Theorem~\ref{minsq.flow} can also be derived from the
network duality in discrete convex analysis (Section 9.6 of
\cite{Murota03}), see Proposition 7.14 in \cite{Frank-Murota.2}.
Analogously to Remark \ref{Sprime} on a slight extension of 
Theorem~\ref{min-maxW}, Theorem~\ref{minsq.flow} can also be easily extended
to the case when $A'$ is a specified subset of edges of $D$, and we
are interested in a non-negative integer-valued $m$-flow $z$ for which
$\sum [z(a)\sp{2}:  a\in A']$ is minimum.  
\hfill $\bullet$ 
\end{remark}

\begin{remark} \rm 
We worked out the details of min-max formulas concerning
the minimum square-sum of a non-negative $m$-flow.  
It is only a technical matter to derive analogous min-max theorems for the minimum
square-sum of 
\RED{
a 
}
feasible ($=$ $(f,g)$-bounded) integral $m$-flow, 
in particular, a circulation or a maximum $st$-flow.  
Our general framework also permits the derivation of a min-max formula 
for the minimum square-sum of feasible integral tension 
\REV{
($=$ potential-difference), 
}
even in the case when not only the 
\REV{
potential-difference
}%
but the potential itself is required to meet upper and lower bounds.
\hfill $\bullet$ 
\end{remark}

\subsection{Minimizing the weighted square-sum}
\label{SCwtdsqsum}

Technically slightly more complicated, but the same approach works for
the weighted square-sum problem.  
Let $a$ be a positive integer and consider the discrete convex function
\begin{equation}  
 \varphi (k):  = a k\sp{2} \quad (k\in {\bf Z}).  
\label{(akk)} 
\end{equation}

\begin{proposition} [\cite{Frank-Murota.2}] 
\label{PRakk.conj}
The discrete conjugate function $\varphi \sp{\bullet}$ 
of $\varphi$ defined in \eqref{(akk)} is given for integers $\ell$ by the following:
\begin{equation}
\label{(akk.conj)} 
 \varphi \sp{\bullet}(\ell) 
 = \llfloor  \frac{\ell + a}{2a}\rrfloor  \ 
 \left( \ell -  a \llfloor  \frac{\ell + a}{2a}\rrfloor  \right). 
\end{equation}
\end{proposition}

\Proof 
The right derivative $\varphi '$ of $\varphi$ 
is given by
$\varphi '(k):= \varphi (k+1)-\varphi (k) = a(k+1)\sp{2} - ak\sp{2} = a(2k+1)$. 
The maximum of $k \ell - \varphi(k)$ is attained by $k$ such that
$\varphi '(k-1) \leq \ell \leq  \varphi '(k)$,
that is,
$ a(2k-1) \leq \ell \leq a(2k+1)$.
Since this is equivalent to 
$(\ell - a )/ (2a) \leq k \leq  (\ell + a)/ (2a)$,
we may take
$k\sp{*}  = \lfloor  (\ell + a)/ (2a ) \rfloor$.
Then 
$\varphi \sp{\bullet}(\ell) = k\sp{*} \ell - \varphi(k\sp{*})$,
which is equal to the right-hand side of \eqref{(akk.conj)}. 
\finbox
\medskip

Theorem~\ref{conjmm} can be written in the following more specific form.

\begin{theorem} \label{min-max-cW} 
Let $R:=\{x:  Qx\geq p\}\subseteq {\bf R}\sp{S}$ 
be a box-TDI polyhedron where $Q$ is an integral matrix and $p$ is an integral vector.  
Let $c$ be a positive integral vector in ${\bf Z}\sp{S}$.  
Then 
\begin{align} 
& \min \{ \sum_{s\in S} c(s)z(s)\sp{2} :  z\in \odotZ{R} \} 
\nonumber \\
&=  
\max \{ yp - \sum_{s\in S} 
 \llfloor  \frac{w(s) + c(s)}{2c(s)}\rrfloor  \
 \left(w(s) - c(s)\llfloor  \frac{w(s) + c(s)}{2c(s)} \rrfloor  
 \right), 
\ w=yQ:  \ y\geq 0 \ \ \hbox{\rm integral}\}.  
\label{(dualmaxcW)} 
\end{align}
\noindent 
Moreover, an integral element $z\sp{*}\in \odotZ{R}$ 
is a minimizer of 
\REV{
\eqref{(dualmaxcW)}
}
if and only if there exists a non-negative integral vector 
$y\sp{*}$ (whose components correspond to the rows of $Q$) 
for which the following optimality criteria hold:
\begin{equation} 
 y\sp{*}(Qz\sp{*}-p)=0, 
\label{(optcrit1.ckk} 
\end{equation}
\begin{equation}  
2c(s)z\sp{*}(s)-1 \leq w\sp{*}(s) \leq 2c(s)z\sp{*}(s)+1 \quad \hbox{\rm for each $s\in S$},
\label{(optcrit2.ckk)} 
\end{equation}
where $w\sp{*}:=y\sp{*}Q$.  
The optimal (integral) dual solution $y\sp{*}$ can be chosen in
such a way that the number of its positive components is at most
$2\vert S\vert$.
\finbox 
\end{theorem}

In Theorems \ref{min-maxW}, \ref{minsq.flow}, and \ref{min-max-cW}, 
we derived min-max formulas
concerning the minimum (weighted) square-sum, and these formulas may
be considered more `standard' from a combinatorial optimization
point of view in the sense that they use neither the notion of
discrete conjugate nor the concept of $\Phi$-compatibility:  they look
like classic combinatorial min-max theorems such as the ones of
Egerv\'ary or Tutte--Berge formula.  
That was made possible by a general min-max
formula relying on the concept of discrete conjugate and by the fact
that in the special case of square-sum we could write up the explicit
form of the discrete conjugate.  This approach shows that the general
min-max formula (Theorem \ref{min-max-gen}) 
can be transformed into a `standard'  one 
whenever one is able to write up explicitly the
discrete conjugate of the separable discrete convex function in question.  
Such a min-max theorem is interesting not only from 
an aesthetic point of view but it is a promising starting point to
develop (purely combinatorial) strongly polynomial algorithms.

This is the reason why it is important to develop a kind of calculus
for concrete discrete conjugates.  
For example, what would 
\REV{
Theorem~\ref{min-maxW} (say) look like 
}
if we were interested in the minimum of the function 
$\Phi$ given by $\Phi(z):=c_{1}z + z\sp 2$, 
or more generally,
$\Phi(z):=c_{1}z + \sum [c_{2}(s)z(s)\sp 2:  s\in S]$ 
(to extend Theorem \ref{min-max-cW}),
where $c_{1}$ and 
$c_{2}\geq 0$ are integral vectors?
Or, what is the conjugate of a function $\Phi$ defined by 
$\Phi(z):  = (z-z_0)\sp{2}$ 
where $z_0$ is a given integral vector?  
In Appendix, we have collected some results of this type.

\section{Inverse combinatorial optimization} 
\label{inverz}

Given a linear 
\REV{
weight- or cost-function
}
 $w_{0}$, find a cheapest
$st$-path, a spanning tree, spanning arborescence, perfect matching,
common basis of two matroids, etc.  
These are standard and well-solved combinatorial optimization problems.  
In an inverse combinatorial optimization problem, 
beside $w_{0}$, we are given an input object $z_{0}$
(path, tree, matching) and the objective is to modify $w_{0}$ as little
as possible so that the input object $z_{0}$ becomes a cheapest one with
respect to the new cost-function $w$.  
If $w_{0}$ is integer-valued, one may require 
that the modified $w$ should also be integer-valued, and
in this section we concentrate exclusively on this case.  
There may be various ways to measure the deviation of $w$ from $w_{0}$.  
For example, in 
\REV{
$l_{1}$-norm
}
the deviation is defined by 
$\sum [\vert w(s) - w_{0}(s)\vert :  s\in S]$.  
\REV{
One may consider weighted versions as well, 
when, for example, the deviation is defined by 
$\sum [c_{1}(s) (w_{0}(s)-w(s)):  w_{0}(s)>w(s)] 
  + \sum [c_{2}(s) (w(s)-w_{0}(s)):  w(s)>w_{0}(s)]$, 
where $c_{1}(s)$ and $c_{2}(s)$ are non-negative integers.  
The $l_{2}$-norm, possibly weighted, is also a natural choice
for measuring the deviation.
}
Even more, imposing lower and upper bounds for the 
\REV{
desired
}
$w$ is also a natural requirement, or, instead of a single input $z_{0}$, 
we may have an input set $\{ z_{1},\dots ,z_{k} \}$ of
solutions and want to find $w$ in such a way that each $z_{i}$ 
is a $w$-minimizer and the deviation of $w$ from $w_{0}$ is minimum.  
Several further versions of 
inverse combinatorial optimization problems have been investigated.  
A relatively early survey paper \cite{Heuberger02} is due to Heuberger, 
\REV{
while the work of Demange and Monnot \cite{DMinv2014}
includes recent developments.
}
Note that Corollary \ref{inverz2} may be viewed
as a solution to a feasibility-type inverse optimization problem.

In this section, we show that the framework in previous sections for
minimizing separable discrete convex functions over a discrete box-TDI set 
covers and even extends an essential part of inverse combinatorial
optimization problems.
Here we concentrate exclusively on the theoretical background and establish
a min-max theorem for the minimum deviation, where the deviation is
measured by an arbitrary separable discrete convex function.  
Our hope is that this theoretical background will provide a good service in
developing efficient algorithms to compute the 
\REV{
desired
}
optimal modification of the input cost-function $w_{0}$.
\REV{
We remark  that in a recent paper  by
Frank and Hajdu \cite{Frank+Hajdu}
 (independently of the present work),
a min-max formula and a simple algorithm have been developed 
for the inverse arborescence problem.  
}%

\subsection{A general framework for inverse problems}
\label{SCinvGenFram}

Let $Qx\geq p$ be a box-TDI system and $R=\{x:Qx\geq p\}$ an integral polyhedron.  
As before, the columns of $Q$ are associated with the elements of ground-set $S$.  
Let $z_{0}\in \odotZ{R}$ be a specified element.  

\REV{
Let $\Phi(w)$ be 
a separable discrete convex function 
\BLU{
on cost-vectors $w$
}
defined as 
$\Phi(w)= \sum_{s\in S} \varphi_{s}(w(s))$
with integer-valued discrete convex functions
$\varphi_{s}$ for $s\in S$.
Let $\ell:  S\rightarrow {\bf Z} \cup \{-\infty \}$ 
and $u:  S\rightarrow {\bf Z} \cup \{+\infty \}$ 
be integral vectors on $S$ 
with $\ell\leq u$,
which represent an interval of admissible cost-vector $w$.
}

The 
\REV{
{\bf inverse separable discrete convex problem}
}
 seeks 
for an integer-valued cost-vector (objective function) $w$ on $S$ 
for which $z_{0}$ is a $w$-minimizer of $R$
 (that is, $w z_{0} \leq w x$ for every $x\in R$), 
$\ell\leq w \leq u$ and $\Phi(w)$ is minimum.
In Corollary~\ref{inverz2}, 
we provided a necessary and sufficient condition 
for the existence of a cost-function $w$ on $S$ for which 
$\ell\leq w \leq u$ and $z_{0}$ is a $w$-minimizer of $\odotZ{R}$.  
Observe that the bounding vectors $\ell$ and $u$ can easily be built 
into $\Phi$ by changing $\varphi_{s}(k)$ to $+\infty$ 
whenever $k>u(s)$ or $k<\ell(s)$ \ ($s\in S$), 
and hence we do not have to work explicitly
with the bounding vectors $\ell$ and $u$.

Our main goal is to characterize those (linear) cost-functions $w$ for
which the input $z_{0}$ is a $w$-minimizer over $R$ and $\Phi(w)$ is
minimum.  We emphasize that $\Phi$ is integer-valued (along with the
bounds $\ell$ and $u$ that can be built into $\Phi$) and 
\REV{
require
}
that the 
\REV{
desired
}
optimal cost-function $w$ is also integer-valued.

In the standard inverse combinatorial optimization problem, as
indicated above, the goal is to modify a starting cost-function $w_{0}$
as little as possible in $l_{1}$-norm so that the input $z_{0}\in R$ is a
$w$-minimizer, where $w$ is the new cost-function.  
For $s\in S$, let $\varphi_{s}(k):= \vert w_{0}(s)-k \vert$.
Then a solution
to the general inverse problem (which minimizes $\Phi$) will provide the 
\REV{
desired
}
solution $w$ for the standard problem.  
With an analogous approach, the general inverse problems can also be built into our
framework of minimizing $\Phi$ over a discrete box-TDI set.  
As a result, the deviation of $w$ from the starting $w_{0}$ may be measured in other norms.  
Moreover, instead of a single initial cost-function $w_{0}$, 
we may specify an interval $[\ell_{0}(s),u_{0}(s)]$ for each $s\in S$ 
and strive to minimize the total deviation of the 
\REV{
desired
}
$w$ from the box defined by these intervals.

\subsection{Preparation}
\label{SCinvPrepa}

In order to embed the general inverse problem into the framework of
discrete box-TDI sets and apply then the min-max results of Section
\ref{main-min-max}, we overview some further properties of box-TDI systems and polyhedra.  
\REV{
Let $C:=\{x:  Kx\geq 0\}$, which is a cone described by an inequality system $Kx\geq 0$,
and let $C\sp{*}$ denote the {\bf dual cone} of $C$, that is, 
$C\sp{*} := \{ w:  w=yK, y\geq 0\}$.  
The {\bf polar cone} of $C$ is $-C\sp{*}$.
}

\begin{proposition} [Chervet, Grappe, Robert \cite{CGR}, Lemma 6]   \label{dualcone} 
A cone is box-TDI if and only if its dual cone is box-TDI.  
\finbox 
\end{proposition}

\begin{proposition} [\cite{CGR}, \REV{Lemma 6}] \label{cone-dilat}
 An integer cone is box-TDI if and only if it is box-integer.  
\finbox 
\end{proposition}

By specializing Theorem~\ref{min-max-gen} to the case of box-TDI cones
\REV{
and using Proposition \ref{cone-dilat},
}
we obtain the following.

\begin{theorem} \label{min-max-cone} 
Let $C$ be a 
\REV{
box-integer
}
cone and let $C\sp{*}$ denote its dual cone.  
Let $\Phi$ be an integer-valued separable discrete convex function on ${\bf Z} \sp{S}$.  
Then 
\begin{align}
 & \min \{ \Phi (z) :  z\in \odotZ{C} \}  
\nonumber \\ &
= \max \{ \Phi(z) - wz :  \ z\in \odotZ{C} , \ w\in \odotZ{C\sp{*}}, 
\ \
\REV{
( z,w ) 
}
\ \ \hbox{\rm $\Phi$-fitting} \}  
\nonumber \\ &
=\max \{ - \Phi\sp{\bullet}(w):  w\in \odotZ{C\sp{*}} \}.
\label{(minmax-cone)} 
\end{align}
\noindent 
An element $z\sp{*}\in \odotZ{C}$ is a $\Phi$-minimizer if and
only if there exists a $w\sp{*}\in \odotZ{C\sp{*}}$ for which $w\sp{*}z\sp{*}=0$ and
\begin{equation}  
 \Phi' (z\sp{*}- \underline{1} ) \leq w\sp{*} \leq \Phi' (z\sp{*}).
\label{(optcrit2x} 
\end{equation} 
\finbox
\end{theorem}

Note that we defined cone $C$ 
\REV{
in terms of its polyhedral description
}
but in the present
formulation we did not make use of this description of $C$.
Therefore, by relying on Proposition~\ref{dualcone}, 
Theorem~\ref{min-max-cone} can be applied to the dual cone $C\sp{*}$ of $C$.

\begin{theorem} \label{min-max-dualcone} 
Let $C$ be a 
\REV{
box-integer
}
cone
and let $C\sp{*}$ denote its dual cone.  
Let $\Phi$ be an
integer-valued separable discrete convex function on ${\bf Z} \sp{S}$.
Then
\begin{align}
& \min \{ \Phi (w) :  w\in \odotZ{C\sp{*}} \} 
\nonumber
\\ & =
\max \{ \Phi(w) - zw :  \ w\in \odotZ{C\sp{*}} , 
\ z\in \odotZ{C}, \ \ 
\REV{
( w, z ) 
}
\ \ \hbox{\rm $\Phi$-fitting} \}  
\nonumber
\\ & =
\max \{ - \Phi\sp{\bullet}(z):  z\in \odotZ{C} \}.  
\label{(minmax-dualcone)} 
\end{align}
\noindent 
An element $w\sp{*}\in \odotZ{C\sp{*}}$ is a $\Phi$-minimizer
if and only if there exists a 
$z\sp{*}\in \odotZ{C}$ for which $w\sp{*}z\sp{*}=0$ and
\begin{equation}
  \Phi' (w\sp{*}-\underline{1}) \leq z\sp{*} \leq \Phi' (w\sp{*}).
\label{(optcrit2xx} 
\end{equation}
\finbox 
\end{theorem} 

\REV{
The following facts will be used in Section \ref{SCinvMinMax}.
}%

\begin{proposition} [\cite{CGR}] \label{boxface} 
Let $x_{1}$ be a solution to a box-TDI system $Qx\geq p$.  
Let $Q_{1}x \geq p_{1}$ denote the subsystem of $Qx\geq p$ 
consisting of those inequalities which are met by $x_{1}$ with equality.  
Then the system $Q_{1}x \geq p_{1}$ is box-TDI.  
\finbox 
\end{proposition}

\begin{remark} \rm \label{tangent}  
\RED{
The polyhedron $C_{1}:=\{x:Q_{1}x\geq p_{1}\}$
is called in \cite{CGR} a ``tangent cone" of $R$.  
Since $C_{1}$ is actually not a cone 
(in the standard meaning of a cone) 
but the translation of cone
$C:=\{x:Q_{1}x\geq 0\}$  by vector $x_{1}$,
we use in this remark the term ``tangent-cone."  
Now Proposition~\ref{boxface} is equivalent to stating that $C_{1}$ is box-TDI, 
which was formulated in
Lemma 5 of \cite{CGR} for minimal ``tangent-cones" of $R$.  
But the proof of Lemma 5 in \cite{CGR} 
works word for word for arbitrary ``tangent-cones" of $R$.
}
\hfill $\bullet$ 
\end{remark}

\begin{proposition} \label{box-cone2} 
Let 
\RED{
$Q$, $p$, $p_{1}$, $x_{1}$ 
}
be the same as in Proposition {\rm \ref{boxface}}.  
Then the cone $C=\{x:  Q_{1}x \geq 0 \}$ is box-TDI.  
\end{proposition}

\Proof 
As mentioned in Remark \ref{tangent},  $C$ is a translation of
$C_{1}=\{x:  Q_{1}x \geq p_{1}\}$.  
By Proposition~\ref{boxface}, 
$C_{1}$ is box-TDI and hence 
\REV{
Proposition~\ref{eltol}
}
implies that $C$ is also box-TDI.  
\finbox.


\subsection{Min-max theorem for the general inverse problem}
\label{SCinvMinMax}

Recall that an ordered pair 
\REV{
$( w, z )$ 
}
of vectors from ${\bf Z}\sp{S}$
\REV{
is
}
called $\Phi$-fitting if 
$\Phi' (w-\underline{1}) \leq z \leq \Phi' (w)$.
Note that we introduced this notion in Section \ref{main-min-max} 
for 
$( z,w )$
 but we use it here for 
$( w,z )$. 
The following result provides a min-max formula for the minimum 
in the inverse separable discrete convex optimization problem 
in which we want to determine the minimum of $\Phi(w)$ 
over those integer-valued linear objective functions $w$ 
for which the input vector $z_{0}\in \odotZ{R}$ minimizes
$wx$ over $R$, that is, $wz_{0} \leq wx$ for each $x\in R$.  
Note that the total dual integrality of the system $Qx\geq p$ implies that
$z_{0}\in \odotZ{R}$ minimizes $wx$ over $R$ 
if and only if 
$z_{0}$ minimizes $wx$ over $\odotZ{R}$.  
We also remark that the duality theorem of linear programming 
implies that $z_{0}$ minimizes $wx$ over $R$ 
if and only if 
$w$ belongs to the cone 
\REV{
$C_{0}\sp{*}$
}
 generated by those rows ${}_{i}q$ 
of $Q$ for which ${}_{i}q z_{0}= p(i)$.

\begin{theorem} \label{gen-inv} 
Let $Qx\geq p$ be a box-TDI system defining
the integral box-TDI polyhedron $R=\{x:  Qx\geq p\}$, 
and let $\Phi$ be an integer-valued separable discrete convex function on ${\bf Z}\sp{S}$.  
Let $z_{0}\in \odotZ{R}$ and let $Q_{0}x\geq p_{0}$ be the subsystem
of $Qx\geq p$ consisting of those inequalities which are met by $z_{0}$ with equalities.  
\REV{
Let $C_{0} :=\{x:  Q_{0}x\geq 0\}$ and let 
$C_{0}\sp{*}:=\{ w: w= y_{0} Q_{0}, y_{0} \geq 0\}$ be the dual cone of $C_{0}$.  
Then 
\begin{align}
& \min \{ \Phi(w): z_{0} \ \hbox{\rm is a $w$-minimizer of}\ \ \odotZ{R}, \ w \
\hbox{\rm integer-valued} \}  
\nonumber
\\ & =
 \max \{ \Phi(w) - zw :  \ w\in \odotZ{C_{0}\sp{*}} , \ z\in \odotZ{C_{0}}, 
 \ \ 
( w, z ) 
\ \ \hbox{\rm $\Phi$-fitting} \} 
\nonumber
\\ & =
 \max \{ - \Phi\sp{\bullet}(z):  z\in \odotZ{C_{0}}\}.  
\label{(minmax-gen-inv)} 
\end{align}
\noindent 
}
An integral cost-function $w\sp{*}$ for which $z_{0}$ is a
$w\sp{*}$-minimizer over $R$ is a $\Phi$-minimizer 
if and only if 
there exists a 
\REV{
$z\sp{*}\in \odotZ{C_{0}}$
}
 for which $w\sp{*}z\sp{*}=0$ 
and the ordered pair 
\REV{
$( w\sp{*},z\sp{*} )$ 
}
is $\Phi$-fitting.  
\end{theorem}

\Proof 
\REV{
As we mentioned before the theorem,
}
$z_{0}$ is a $w\sp{*}$-minimizer element of $R$ 
precisely if 
\REV{
$w\sp{*}\in C_{0}\sp{*}$.  
}
\REV{
By applying Proposition~\ref{box-cone2} to
$C_{0}$ in place of $C$, we obtain that $C_0$ 
}
is box-TDI and hence 
Theorem~\ref{min-max-dualcone} implies the theorem.  
\finbox

\medskip

\REV{
\begin{remark} \rm \label{RMinvwkassump}  
The proof of Theorem~\ref{gen-inv} shows that the role of the assumed box-TDI-ness
of $R$ is only to ensure that $C_{0}$ is box-TDI (equivalently, box-integer).
If $z_{0}$ is specified and fixed,
we can weaken the assumption to box-TDI-ness of the tangent cone of $R$ at $z_{0}$.
\hfill $\bullet$ 
\end{remark}
}

\REV{
\begin{remark} \rm \label{RMinveqn}  
Theorem~\ref{gen-inv} can easily be extended for the case 
when the box-TDI system defining $R$ is given in the form of
$[Q'x\geq p', Q\sp{=} x = p\sp{=} ]$.  
In this case, let $Q'_{0} x\geq p'_{0}$ be the subsystem of $Q'x\geq p'$ 
consisting of those inequalities which are met by $z_{0}$ with equalities.  
Then 
\eqref{(minmax-gen-inv)}
holds for 
$C_{0}:=\{x:  Q'_{0}x\geq 0, Q\sp{=} x=0\}$ and its dual cone
$C_{0}\sp{*}:=\{w:  w=y'_{0}Q'_{0}+ y_{0}\sp{=} Q\sp{=} , y'_{0}\geq 0\}$. 
\hfill $\bullet$ 
\end{remark}
}

\begin{remark} \rm \label{RMinvext}  
A natural extension of the problem is when, 
instead of a single element $z_{0}$, 
we have a subset $Z_{0}:=\{z_{1},z_{2},\dots ,z_{k}\}$ of
elements of $\odotZ{R}$, and the goal is to characterize those
integer-valued weight-functions $w$ for which each $z_{i}\in Z_{0}$ 
is a $w$-minimizer element of $\odotZ{R}$ and $\Phi(w)$ is minimum. 
 (It is allowed that $\odotZ{R}$ may have other $w$-minimizer elements.) 
To treat this case let $R_{k}:=kR$ denote the $k$-dilation of $R$.  
By Proposition~\ref{dilat}, $R_{k}$ is also a box-TDI polyhedron containing
$z_{0}:=z_{1}+\cdots +z_{k}$.  
It is a straightforward observation for a cost-function $w$ that $z_{0}$ 
is a $w$-minimizer element of $R_{k}$
precisely if each $z_{i}$ is a $w$-minimizer of $R$.  
Therefore we can
apply Theorem~\ref{gen-inv} to $k$-dilation $R_{k}$ of $R$ and to
$z_{0}:=z_{1}+\cdots +z_{k}$.
\hfill $\bullet$ 
\end{remark}

\begin{remark} \rm \label{RMappe}  
In Appendix we overview some
special separable discrete convex functions related to (weighted)
$l_{1}$-norm, and calculate their explicit discrete conjugates,
analogously to the way how the conjugate of the (weighted) square-sum
was calculated in Section~\ref{SCspedcfn}. 
By applying Theorem~\ref{gen-inv} to these concrete conjugates, one can obtain min-max
formulas of standard combinatorial optimization type (that is, without
using conjugate) for a great number of inverse problems.  One example
is the inverse matroid intersection problem when there is a specified
upper and lower bound for the 
\REV{
desired
}
cost-function $w$.  
In another version, we want to minimize the deviation of the 
\REV{
desired
}
cost-function from a specified box, rather than from a single point $w_{0}$.
With this framework, one can derive min-max theorems even 
for minimum cost versions of the inverse problems where a linear cost-function is
specified for the deviation of $w$ from $w_{0}$.  Beyond theoretical
advantage, our hope is that this kind of min-max formulas shall
facilitate the development of strongly polynomial algorithms for these
cases.  
\hfill $\bullet$ 
\end{remark}

\section{Appendix: Calculating concrete discrete conjugates}
\label{SCconcrconj}

Recall that a function 
$\varphi :{\bf Z} \rightarrow {\bf Z} \cup \{+\infty \}$ 
is called discrete convex if 
$\varphi (k-1) + \varphi (k+1) \geq 2\varphi (k)$ 
for each $k\in $ dom$(\varphi )$.  
Below we list the discrete conjugate of some concrete univariate functions and
for some elementary operations.  We emphasize that every function is
assumed to be integer-valued.
Naturally, these formulas immediately
extend to separable discrete convex functions.  
The proof of these claims are not difficult and left to the reader.

\begin{claim}[\cite{Murota98a,Murota03}] \label{bicon} 
The discrete conjugate of a discrete convex function $\varphi$ is discrete convex.  
Furthermore,
\begin{equation}
 (\varphi \sp{\bullet})\sp{\bullet} = \varphi . 
\label{(bicon)} 
\end{equation} 
\end{claim}

\begin{claim}[\cite{Murota98a,Murota03}] \label{sum} 
The sum $\varphi =\varphi_{1}+\varphi_{2}$
of two discrete convex functions is discrete convex 
and its discrete conjugate $\varphi \sp{\bullet}$ is given by
\begin{equation} 
\varphi \sp{\bullet}(\ell) = 
  \min \{\varphi_{1}\sp{\bullet}(\ell_{1}) + \varphi_{2}\sp{\bullet}(\ell_{2}):  
         \ell_{1}+\ell_{2} = \ell, \ \ell_{i}\in {\bf Z}  \}.
\label{(sum.1)} 
\end{equation}
\noindent 
In the special case when $\varphi_{2}$ is a linear function
defined by $\varphi_{2}(k):= ck$, where $c$ is an integer, one has
\begin{equation} 
\varphi \sp{\bullet}(\ell)= \varphi_{1}\sp{\bullet}(\ell-c).  
\label{(sum.2)}
\end{equation}
\end{claim}

\begin{claim}  \label{k-k0} 
For an integer $k_{0}$, the discrete conjugate of a
discrete convex function $\varphi_{0}$ defined by 
$\varphi_{0}(k):=\varphi (k-k_{0})$ 
is given by
\begin{equation} 
\varphi_{0}\sp{\bullet}(\ell) = \varphi \sp{\bullet}(\ell) + k_{0}\ell.  
\label{(k-k0)}
\end{equation} 
\end{claim}

The next claim is useful in situations when we want to build specified
lower and upper bounds imposed on the variables 
into the function $\varphi $ to be minimized, 
by making $\varphi $ to be $+\infty $ outside the bounds.

\begin{claim}  \label{Irestr} 
Let $A \leq B $ be integers, 
where $A $ may be $-\infty $ and $B $ may be $+\infty $, 
and let $I:=[A ,B ]_{\bf Z} $ be
the set of integers $k$ with $A \leq k\leq B $.  
Let $\varphi_{I}$ denote the function obtained from a 
discrete convex function $\varphi $ 
by restricting it to $I$ in the following sense:
\begin{equation} 
\varphi_{I}(k) = 
\begin{cases} 
\varphi (k) & \quad \hbox{\rm if}\quad A \leq k\leq B  , 
\cr 
+\infty & \quad \hbox{\rm otherwise} . 
\end{cases}
\label{(Irestr)} 
\end{equation}
\noindent 
Then $\varphi_{I}$ is discrete convex and its discrete conjugate is as follows:
\begin{equation}
\varphi_{I}\sp{\bullet}(\ell) 
= \min \{ \varphi (\ell_{1}) + \max \{A \ell_{2}, B \ell_{2}\}:  
      \ell_{1}+\ell_{2} = \ell, \  \ell_{i} \ \hbox{\rm integer} \}.
\label{(Irestr2)} 
\end{equation} 
\end{claim}

\begin{claim}  \label{k0AB} 
Let $c_{-} \leq c_{+} $ be integers, and let $A \leq k_{0}\leq B $ be integers, 
where $A $ may be $-\infty $ and $B $ may be $+\infty $. 
Let $\varphi $ be defined by
\begin{equation} 
\varphi (k):  = \begin{cases} 
c_{-} (k-k_{0}) & \quad \hbox{\rm if}\quad A \leq k\leq k_{0} ,
\cr 
c_{+} (k-k_{0}) & \quad \hbox{\rm if}\quad k_{0}\leq k\leq B  ,
\cr 
+\infty & \quad \hbox{\rm otherwise}.  
\end{cases} 
\label{(k0AB.1)} 
\end{equation}
\noindent 
Then the discrete conjugate of $\varphi $ is as follows:
\begin{equation} 
\varphi \sp{\bullet}(\ell) = \begin{cases}
A \ell - c_{-} (A -k_{0}) & \quad \hbox{\rm if}\quad \ell < c_{-} ,
\cr 
k_{0}\ell & \quad \hbox{\rm if}\quad c_{-} \leq \ell \leq c_{+}  ,
\cr
B \ell - c_{+} (B - k_{0}) & \quad \hbox{\rm if} \quad \ell > c_{+} .
\end{cases} 
\label{(k0AB.2)} 
\end{equation}
\noindent 
When $A =-\infty $, for case $\ell < c_{-} $ one has 
$\varphi\sp{\bullet}(\ell) = A \ell - c_{-} (A -k_{0}) 
= A  (\ell - c_{-} ) + c_{-} \, k_{0} = +\infty $,
and, analogously, when $B =+\infty $, for case $\ell > c_{+}$ 
one has $\varphi \sp{\bullet}(\ell) = +\infty $. 
\end{claim}

\begin{claim} \label{abAB} 
Let $c_{-} \leq 0\leq c_{+} $ be integers.  
Let $A \leq a <b \leq B $ 
be integers where $A $, $a$ may be $-\infty $ 
and $B $, $b$ may be $+\infty $. 
Let $\varphi $ be defined by:
\begin{equation} 
\varphi (k):  = \begin{cases}
 0 & \quad \hbox{\rm if}\quad a \leq k\leq b ,
\cr 
c_{-} (k- a) & \quad \hbox{\rm if}\quad A \leq k < a , 
\cr 
c_{+} (k- b) & \quad \hbox{\rm if}\quad b <k\leq B  ,
\cr
+\infty & \quad \hbox{\rm otherwise}.
\end{cases} 
\label{(abAB.1)} 
\end{equation}
\noindent 
Then the discrete conjugate of $\varphi $ is as follows:
\begin{equation} 
\varphi \sp{\bullet}(\ell) = \begin{cases} 
A \ell - c_{-} (A - a) & \quad \hbox{\rm if}\quad \ell < c_{-} ,
\cr 
a \ell & \quad \hbox{\rm if}\quad c_{-} \leq \ell < 0 ,
\cr
 0 & \quad \hbox{\rm if}\quad \ell = 0 ,
\cr 
b \ell & \quad \hbox{\rm if}\quad 0< \ell \leq c_{+} , 
\cr 
B \ell - c_{+} (B - b) & \quad \hbox{\rm if}\quad \ell > c_{+} .  
\end{cases} 
\label{(abAB.2)} 
\end{equation}
\noindent 
When $A =-\infty $, for case $\ell < c_{-} $ one has 
$\varphi\sp{\bullet}(\ell) = A \ell - c_{-} (A - a) 
= A  (\ell - c_{-} ) +c_{-} a = +\infty $, 
and, analogously, 
when $B =+\infty $, for case $\ell > c_{+} $ 
one has $\varphi \sp{\bullet}(\ell) = +\infty $. 
\end{claim}


\paragraph{Acknowledgement} \ 
The authors are grateful to R. Grappe for his
indispensable and profound help concerning fundamental properties of
box-TDI polyhedra.  
The detailed comments from the referees were helpful to improve the paper.
The research was partially supported by the
National Research, Development and Innovation Fund of Hungary (FK-18)- No.~NKFI-128673, 
and by JSPS KAKENHI Grant Number JP20K11697.




\begin{thebibliography}{999}
\setlength{\itemsep}{0pt}





\bibitem{BL}
J.M. Borwein and A.S. Lewis, 
Convex Analysis and Nonlinear Optimization, Theory and Examples, 
(Second Edition) 2005, Canadian
Mathematical Society. CMS Books in Mathematics.



\bibitem{Cameron89}
K. Cameron, 
{\it A min-max relation for the partial $q$-colourings of a graph, Part II:
Box perfection}, Discrete Mathematics, 74 (1989), 15--27.


\bibitem{CGR}
 P. Chervet, R. Grappe, L.-H.  Robert, 
{\it Box-total dual integrality, box-integrality, and equimodular matrices}, 
Mathematical Programming, Ser.  A, published online:  20 May 2020.
https://doi.org/10.1007/s10107-020-01514-0



\bibitem{Cook83}
 W.J. Cook, 
{\it Operations that preserve total dual integrality}, 
Operations Research Letters, 2 (1983) 31--35.

\bibitem{Cook86}
 W. Cook, 
{\it On box totally dual integral polyhedra}, 
Mathematical Programming,  34 (1986) 48--61.

\bibitem{CFS}
 W.J.  Cook, J. Fonlupt, and A. Schrijver, 
{\it An integer analogue of Carath{\'e}odory's theorem}, 
J. Combinatorial Theory, Ser B. 40 (1986) 63--70.



\bibitem{DMinv2014}
M. Demange and J. Monnot, 
{\it An introduction to inverse combinatorial problems}, 
Chapter~17 in:  
Paradigms of Combinatorial Optimization: Problems and New Approaches, 
second edition
(V.Th.~Paschos, ed.), 
ISTE LTd and John Wiley and Sons, 
2014. pp.~547--586.




\bibitem{Edmonds70}
 J. Edmonds, 
{\it Submodular functions, matroids, and certain polyhedra}, 
in:   Combinatorial Structures and their Applications
 (R.  Guy, H. Hanani, N. Sauer, and J. Sch\"onheim, eds.), 
Gordon and Breach, New York (1970) pp. 69--87.


\bibitem{Edmonds71}
 J. Edmonds, 
{\it Matroids and the greedy algorithm}, 
Mathematical Programming, 1 (1971) 127--136.


\bibitem{Edmonds-Giles}
 J. Edmonds and R. Giles, 
{\it A min-max relation for submodular functions on graphs}, 
Annals of Discrete Mathematics, 1 (1977), 185--204.



\bibitem{Edmonds-Giles84}
 J. Edmonds and R. Giles, 
{\it Total dual integrality of linear inequality systems}, 
in:  Progress in Combinatorial Optimization
(ed.  W. R. Pulleyblank) Academic Press
(1984) 117--129.




\bibitem{Frank-book}
 A. Frank, 
Connections in Combinatorial Optimization, \ Oxford University Press, 
2011 
(ISBN 978-0-19-920527-1),
Oxford Lecture Series in Mathematics and its Applications, 38.



\bibitem{Frank+Hajdu}
A. Frank and G. Hajdu,
{\it A simple algorithm and min-max formula for the inverse arborescence problem}, 
to appear in Discrete Applied Mathematics.





\bibitem{Frank-Murota.2}
A. Frank  and K. Murota,
{\it Discrete decreasing minimization, Part II:
Views from discrete convex analysis}, 
arXiv: 1808.08477v4 30, June 2020.


\bibitem{FM19partA} 
A. Frank and  K.  Murota,
{\it Decreasing minimization on M-convex sets},
arXiv: 2007.09616, July 2020.



\bibitem{Fujishige80} 
S. Fujishige, 
{\it Lexicographically optimal base of a polymatroid with respect to a  weight vector},
 Mathematics of Operations Research, 5 (1980) 186--196.




\bibitem{Groenevelt91}
 H. Groenevelt, 
{\it Two algorithms for maximizing a separable concave function 
over a polymatroid feasible region}, 
European J. of Operational Research, 54 (1991) 227--236.



\bibitem{Heuberger02}
 C. Heuberger, 
{\it Inverse combinatorial optimization:  
A survey on problems, methods, and results}, 
Journal of Combinatorial Optimization, 8 (2004) 329--361.


\bibitem{HL01}
J.-B. Hiriart-Urruty and C. Lemar{\'e}chal,
Fundamentals of Convex Analysis, 
Springer, Berlin, 2001.





\bibitem{KOS} V. Kaibel, S. Onn, P. Sarrabezolles, 
{\it The unimodular intersection problem}, 
Operations Research Letters,  43 (2015), 502--504.





\bibitem{McC06handbk} 
S. T. McCormick,
{\it Submodular function minimization}, 
in: K. Aardal, G. Nemhauser, and R. Weismantel (Eds.),
Handbook on Discrete Optimization,
Elsevier Science Publishers, Berlin, 2006, Chapter 7, pp.321--391.



\bibitem{MPSV20}
S.T.  McCormick, B. Peis, R. Scheidweiler, and F. Valentin, 
{\it A polynomial time algorithm for solving the closest
vector problem in zonotopal lattices}, 
arXiv: 2004.07574v1 
April 2020.


\bibitem{Murota98a} 
K. Murota, {\it Discrete convex analysis}, 
Mathematical Programming, 83 (1998) 313--371.


\bibitem{Murota03}
 K. Murota, Discrete Convex Analysis,
SIAM, Philadelphia, 2003.










\bibitem{Roc70}
R.T. Rockafellar, Convex Analysis,
Princeton University Press, Princeton, 1970.




\bibitem{Schrijverbook0}
 A. Schrijver, Theory of Linear and Integer Programming,
 Wiley, Chichester, 1986.

\bibitem{Schrijverbook}
 A. Schrijver, 
 Combinatorial Optimization:  Polyhedra and Efficiency, Springer, Heidelberg, 2003.



\bibitem{Sebo90} 
A. Seb{\H o}, 
{\it Hilbert bases, Caratheodory's theorem and combinatorial optimization}, 
in:  R. Kannan, W.R. Pulleyblank (Eds.), 
Integer Programming and Combinatorial Optimization, 
University of Waterloo Press, Waterloo, Canada, 1990, pp.431--455.








\end{thebibliography}
\end{document}